\documentclass[10pt]{article}

\usepackage{a4wide}
\usepackage{amssymb}
\usepackage{amsfonts}
\usepackage{amsmath}
\input xy
\xyoption{arrow} \xyoption{matrix}

\date{}

\newtheorem{proposition}{Proposition}[section]
\newtheorem{theorem}[proposition]{Theorem}
\newtheorem{lemma}[proposition]{Lemma}

\newtheorem{corollary}[proposition]{Corollary}

\def\GK{{\rm  GK}\,}

\def\der{\partial }

\def\nFM0{{\nu }_{F,M_0}}
\def\nFN0{{\nu }_{F,N_0}}
\def\nGN0{{\nu }_{G,N_0}}

\def\N0{ {\bf N}_0 }

\def\t{\otimes}
\def\g{\gamma}
\def\v{\varphi}
\def\ra{\rightarrow}

\def\lra{\leftrightarrow}
\def\Xpm{X^{\pm }}

\def\s{\sigma}
\def\Z{\mathbb{Z}}

\def\l1{{\lambda}_1}

\def\a{\alpha}
\def\a0{ {\alpha }_0}
\def\a1{ {\alpha }_1}

\def\l{\lambda}

\def\nFGM0{{\nu }_{F,G,M_0}}

\def\nFN0{{\nu}_{F,N_0}}

\def\sm{{\sigma}^m}

\def\sm1{{\sigma}^{-1}}

\def\smtp1{{\sigma}^{-t+1}}

\def\S1{S^{-1}}

\def\Xpm1{X^{\pm 1}_1}

\def\sPM1{{\sigma }^{\pm 1}}
\def\sMP1{{\sigma }^{\mp 1 }}

\def\d{\delta}

\def\fd{{\rm fd}_A}

\def\di{{\rm d.ind}}

\def\L{\Lambda}

\def\CA{{\cal A}}

\def\Ytm1{Y^{t-1}}
\def\Yim1{Y^{i-1}}

\def\CL{{\cal L}}
\def\CM{{\cal M}}

\def\CH{{\cal H}}
\def\ass{{\rm ass}}

\def\ad{{\rm ad }}
\def\dim{{\rm dim }}

\def\ker{ {\rm ker } }

\def\CJ{ {\cal J}}
\def\gr{ {\rm gr} }

\def\SL2Z{ {\rm SL}_2({\bf Z}) }

\def\CL{{\cal L}}

\def\Gp1{ G^{1 , 1 } }
\def\P11{ P^{-1 , 1 } }
\def\Pp1{ P^{1 , 1 } }

\def\nCLsr{{}^\nu\kern-2pt {\cal L}^{\sigma , \rho  }}
\def\nP{{}^\nu \kern-2pt P}
\def\nL{{}^\nu\kern-2pt L}
\def\nLL{{}^\nu\kern-2pt \Lambda}
\def\nPsr{{}^\nu\kern-2pt P^{\sigma , \rho  }}
\def\nLsr{{}^\nu\kern-2pt L^{\sigma , \rho  }}
\def\nuCL{{}^\nu\kern-2pt  {\cal L}}
\def\nCLsr{{}^\nu\kern-2pt {\cal L}^{\sigma , \rho  }}
\def\nCL1m{{}^\nu\kern-2pt {\cal L}^{-1 , 1  }}
\def\x1nu{x^\frac{1}{\nu}}
\def\xm1nu{x^{-\frac{1}{\nu}}}

\def\ra{\rightarrow }

\def\CB{{\cal B}}

\def\CI{{\cal I}}

\def\CT{{\cal T}}

\def\CC{ {\cal C}}

\def\CH{ {\cal H}}
\def\CP{ {\cal P}}

\def\nAM0{{\nu }_{{\cal A},M_0}}
\def\nAN0{{\nu }_{{\cal A},N_0}}

\def\End{ {\rm End }}

\def\CJ{ {\cal J }}

\def\CP{ {\cal P }}
\def\det{ {\rm det }}
\def\ad{ {\rm ad }}

\def\ga{\mathfrak{a}}
\def\gb{\mathfrak{b}}
\def\gc{\mathfrak{c}}
\def\gd{\mathfrak{d}}

\def\gp{\mathfrak{p}}
\def\gq{\mathfrak{q}}
\def\gr{\mathfrak{r}}

\def\GL{{\rm GL}}
\def\SL{{\rm SL}}

\def\Spec{{\rm Spec}}

\def\di!{\frac{\der^i}{i!}}
\def\dik!{\frac{\der^k_i}{k!}}

\def\id{{\rm id}}

\def\gldim{{\rm gldim}}

\def\N{\mathbb{N}}

\def\0{\overline{0}}
\def\1{\overline{1}}

\def\Ln1{\L_{n,\overline{1}}}

\def\oa{\overline{a}}

\def\a1{a_{\overline{1}}}

\def\S{\Sigma}

\def\vn1{\overrightarrow{n-1}}

\def\im{{\rm im}}

\def\mA{\mathbb{A}}
\def\mD{\mathbb{D}}

\def\gf{\mathfrak{f}}

\def\Sub{{\rm Sub}}

\def\Min{{\rm Min}}

\def\mS{\mathbb{S}}

\def\mI{\mathbb{I}}

\def\lann{{\rm l.ann}}
\def\rann{{\rm r.ann}}
\def\Cen{{\rm Cen}}
\def\clKdim{{\rm cl.Kdim}}
\def\pd{{\rm pd}}

\def\lgldim{{\rm l.gldim}}
\def\rgldim{{\rm r.gldim}}

\def\fd{{\rm fd}}
\def\bM{\overline{M}}

\def\mF{\mathbb{F}}

\def\Id{\rm Id}
\def\wdim{{\rm wdim}}

\begin{document}

\author{V. V. \  Bavula 
}

\title{The algebra of integro-differential operators on a
polynomial algebra}

\maketitle

\begin{abstract}
We prove that the   algebra $\mI_n:=K\langle x_1, \ldots , x_n,
\frac{\der}{\der x_1}, \ldots ,\frac{\der}{\der x_n},  \int_1,
\ldots , \int_n\rangle $ of   integro-differential operators  on a
polynomial algebra is a prime, central, catenary, self-dual,
non-Noetherian algebra of classical Krull dimension $n$ and of
Gelfand-Kirillov dimension $2n$. Its weak homological dimension is
$n$, and  $n\leq  \gldim (\mI_n)\leq 2n$. All the ideals of
$\mI_n$ are found explicitly, there are only finitely many of them
($\leq 2^{2^n}$), they commute ($\ga \gb = \gb\ga$) and are
idempotent ideals ($\ga^2= \ga$). The number of ideals of $\mI_n$
is equal to  the {\em Dedekind number} $\gd_n$.
  An analogue of Hilbert's Syzygy
Theorem is proved for $\mI_n$. The group of units of the algebra
$\mI_n$ is described (it is a huge group). A canonical form is
found for each integro-differential operators (by proving that the
algebra $\mI_n$ is a generalized Weyl algebra). All the mentioned
results hold for the Jacobian algebra $\mA_n$  (but $\GK (\mA_n)
=3n$, note that $\mI_n\subset \mA_n$). It is proved that the
algebras $\mI_n$ and $\mA_n$ are ideal equivalent.



$${\bf Contents}$$
\begin{enumerate}
\item Introduction. \item Defining relations for the algebra
$\mI_n$. \item Ideals of  the algebra $\mI_n$. \item The
Noetherian factor algebra of the algebra $\mI_n$. \item The group
of units of the algebra $\mI_n$ and its centre.  \item The weak
and the global dimensions of the algebra $\mI_n$. \item The weak
and the global dimensions of the Jacobian algebra $\mA_n$.
\end{enumerate}
\end{abstract}


\section{Introduction}\label{INTRO}
Throughout, ring means an associative ring with $1$; module means
a left module;
 $\N :=\{0, 1, \ldots \}$ is the set of natural numbers; $K$ is a
field of characteristic zero and  $K^*$ is its group of units;
$P_n:= K[x_1, \ldots , x_n]$ is a polynomial algebra over $K$;
$\der_1:=\frac{\der}{\der x_1}, \ldots , \der_n:=\frac{\der}{\der
x_n}$ are the partial derivatives ($K$-linear derivations) of
$P_n$; $\End_K(P_n)$ is the algebra of all $K$-linear maps from
$P_n$ to $P_n$;
the
subalgebra  $A_n:= K \langle x_1, \ldots , x_n , \der_1, \ldots ,
\der_n\rangle$ of $\End_K(P_n)$ is called the $n$'th {\em Weyl}
algebra.

$\noindent $

{\it Definition}, \cite{Bav-Jacalg}: The {\em Jacobian algebra}
$\mA_n$ is the subalgebra of $\End_K(P_n)$ generated by the Weyl
algebra $A_n$ and the elements $H_1^{-1}, \ldots , H_n^{-1}\in
\End_K(P_n)$ where $$H_1:= \der_1x_1, \ldots , H_n:= \der_nx_n.$$

Clearly, $\mA_n =\bigotimes_{i=1}^n \mA_1(i) \simeq \mA_1^{\t n }$
where $\mA_1(i) := K\langle x_i, \der_i , H_i^{-1}\rangle \simeq
\mA_1$. The algebra $\mA_n$ contains all the  integrations
$\int_i: P_n\ra P_n$, $ p\mapsto \int p \, dx_i$, since  $$\int_i=
x_iH_i^{-1}: x^\alpha \mapsto (\alpha_i+1)^{-1}x_ix^\alpha.$$ In
particular, the algebra $\mA_n$ contains the {\em algebra}
$\mI_n:=K\langle x_1, \ldots , x_n$,
 $\der_1, \ldots ,\der_n,  \int_1,
\ldots , \int_n\rangle $ {\em of  polynomial integro-differential
operators}. Note that $\mI_n=\bigotimes_{i=1}^n\mI_1(i)\simeq
\mI_1^{\t n}$ where $\mI_1(i):= K\langle x_i, \der_i,
\int_i\rangle$.

$\noindent $

The paper proceeds as follows. In Section \ref{DRAI}, two sets of
defining relations are given for the algebra $\mI_n$ (Proposition
\ref{a5Oct9}); a canonical form is found for each element of
$\mI_n$ by showing that the algebra $\mI_n$ is a generalized Weyl
algebra (Proposition \ref{a5Oct9}.(2)); the Gelfand-Kirillov
dimension of the algebra $\mI_n$ is $2n$ (Theorem \ref{5Oct9}).

In Section \ref{IOTAI}, a new equivalence relation, the {\em ideal
equivalence}, on the class of algebras is introduced: two algebras
$A$ and $B$ are {\em ideal equivalent} if there exists a bijection
$f$ from the set $\CJ (A)$ of all the ideals of the algebra $A$ to
the set $\CJ (B)$ of all the ideals of the algebra $B$ such that,
for all $\ga, \gb \in \CJ (A)$,
$$ f(\ga +\gb)= f(\ga )+ f(\gb ), \;\; f(\ga\cap  \gb)= f(\ga ) \cap f(\gb ),
\;\; f(\ga \gb) =f(\ga ) f(\gb ).$$ The algebras $\mI_n$ and
$\mA_n$ are ideal equivalent (Theorem \ref{7Oct9}). As a result,
we have for free many results for the ideals of $\mI_n$ using
similar known results for the ideals of $\mA_n$ of
\cite{Bav-Jacalg}. Name just a few:
\begin{itemize}
\item $\mI_n$ is a prime, catenary algebra of classical Krull
dimension $n$, and there is a unique maximal ideal $\ga_n$ of the
algebra $\mI_n$. \item $\ga \gb=\gb\ga $ and $\ga^2= \ga$ for all
$\ga , \gb \in \CJ (\mI_n)$.  \item The lattice $\CJ (\mI_n)$ is
distributive. \item Classifications of all the ideals and the
prime ideals of the algebra $\mI_n$ are given. \item The set $\CJ
(\mI_n)$ is finite. Moreover, $|\CJ (\mI_n)|=\gd_n$ where $\gd_n $
is the Dedekind number, and $2-n+\sum_{i=1}^n2^{n\choose i}\leq
\gd_n \leq 2^{2^n}$. \item $P_n$ is the only (up to isomorphism)
faithful simple $\mI_n$-module.
\end{itemize}

The fact that certain rings of differential operators are catenary
was proved by  Brown, Goodearl and  Lenagan in
\cite{Brown-Goodearl-Lenagan}.

In Section \ref{NFAAI}, it is proved that the factor algebra
$\mI_n/ \ga$ is Noetherian iff the ideal $\ga$ is maximal
(Proposition \ref{A17Oct9}); and $\GK (\mI_n/\ga)=2n$ for all the
ideals $\ga$ of $\mI_n$ distinct from $\mI_n$ (Lemma
\ref{a18Oct9}).

In Section \ref{GUAIC}, for the algebra $\mI_n$ an involution $*$
is introduced such that $\der_i^*=\int_i$, $\int_i^* = \der_i$,
and $H_i^* = H_i$, see (\ref{*invIn}). This means that the algebra
$\mI_n$ is `symmetrical' with respect to derivations and
integrations. $\ga^* = \ga$ for all ideals $\ga$ of the algebra
$\mI_n$ (Lemma \ref{b12Oct9}.(1)). Each ideal of the algebra
$\mI_n$ is an essential left and right submodule of $\mI_n$ (Lemma
\ref{a17Oct9}.(2)). The group $\mI_n^*$ of units of the algebra
$\mI_n$ is described:
$$ \mI_n^* = K^* \times (1+\ga_n)^* \supseteq K^* \times \underbrace{\GL_\infty
(K)\ltimes\cdots \ltimes \GL_\infty (K)}_{2^n-1 \;\; {\rm
times}}$$ and its centre is $K^*$ (Theorem \ref{17Oct9}). For
$n=1$, the group $\mI_n^*$ is found explicitly, $\mI_1^*\simeq
K^*\times \GL_\infty (K)$ (Corollary \ref{b17Oct9}). The centre of
the algebra $\mI_n$ is $K$ (Lemma \ref{c12Oct9}.(2)). It is proved
that, for a $K$-algebra $A$, the algebra $A\t \mI_n$ is prime iff
the algebra $A$ is prime (Corollary \ref{a24Oct9}).

In Section \ref{WGDAI}, we prove that 
 the weak  (w.dim) dimension of the algebra $\mI_n$ is $n$ (Theorem
\ref{18Oct9}). Moreover, $\wdim (\mI_n / \gp )=n$ for all the
prime ideals $\gp \in \Spec (\mI_n)$ (Corollary \ref{a19Oct9}).
Recall that for each Noetherian ring its weak dimension coincides
with its global dimension (in general, this is wrong for
non-Noetherian rings). In 1972, Roos proved that the global
dimension of the Weyl algebra $A_n$ is $n$, \cite{Roos}. This
result was generalized by Chase \cite{Chase-1974} to the ring of
differential operators on a smooth affine variety. Goodearl
obtained formulae for the global dimension of certain rings of
differential operators \cite{Goodearl-gldim-II},
\cite{Goodearl-gldim-III} (see also Levasseur
\cite{Levasseur-injdim}, and van den Bergh
\cite{VandenBergh-1991}). Holland and Stafford found the global
dimension of the ring of differential operators on a rational
projective curve \cite{Holland-Stafford-1992} (see also Smith and
Stafford \cite{Smith-Stafford-1988}).

 Many classical algebras are tensor product of algebras (eg,
 $P_n=P_1^{\t n}$, $A_n= A_1^{\t n}$, $\mA_n= \mA_1^{\t n}$,
 $\mI_n=\mI_1^{\t n}$, etc). In general, it is difficult to find the
 dimension $d(A\t B)$ of the tensor product of two algebras (even to
 answer the question of  when it is finite). In
 \cite{ERZ},  it was pointed out  by Eilenberg, Rosenberg and  Zelinsky that `{\em
the questions concerning the dimension of the tensor product of
two algebras have turned out to be surprisingly difficult.}' An
answer is known if one of the algebras is a polynomial algebra:
$${\bf  Hilbert's \; Syzygy\;  Theorem}:\;\;\;  d(P_n \t B) = d(P_n) +d(B) =
n+d(B), $$ where $d=\wdim, \, \gldim$. In \cite{THM},
\cite{glgwa}, an analogue of Hilbert's Syzygy Theorem was
established for certain generalized Weyl algebras $A$ (eg,
$A=A_n$, the Weyl algebra):
$$ \lgldim (A\t B) = \lgldim (A) +\lgldim (B)$$
for all left Noetherian finitely generated algebras $B$ ($K$ is an
algebraically closed uncountable field of characteristic zero). In
this paper, a similar result is proved for the algebra $\mI_n$ and
for all its prime factor algebras but for the weak dimension
(Theorem \ref{19Oct9}). It is shown that $n\leq \gldim (\mI_n)
\leq 2n$ (Proposition \ref{a29Nov9}).

In Section \ref{WGDJJ}, we prove that 
 the weak
dimension of the Jacobian algebra $\mA_n$ is $n$ (Theorem
\ref{J18Oct9}), and  $\wdim (\mA_n / \gp )=n$ for all the prime
ideals $\gp \in \Spec (\mA_n)$ (Corollary \ref{Ja19Oct9}). An
analogue of
 Hilbert's Syzygy Theorem is proved for the Jacobian  algebra
$\mA_n$ and its prime factor algebras (Theorem \ref{J19Oct9}). It
is shown that $n\leq \gldim (\mA_n) \leq 2n$ (Proposition
\ref{e29Nov9}).

The algebra $\mI_1= A_1\langle \int \rangle$ is an example of the
{\em Rota-Baxter} algebra. The latter appeared in the work of
Baxter \cite{Baxter} and further explored by Rota \cite{Rota-1,
Rota-2}, Cartier \cite{Cartier}, and Atkinson \cite{Atkinson}, and
more recently by many others: Aguiar, Moreira
\cite{Aguiar-Moreira}; Cassidy, Guo, Keigher, Sit, Ebrahimi-Fard
\cite{Cassidy-Guo-Keigher-Sit}, \cite{Ebrahimi-Guo}; Connes,
Kreimer, Marcoli \cite{Connes-Kreimer}, \cite{Connes-Marcoli},
name just a  few.  From the angle of the Rota-Baxter algebras the
algebra $\mI_1$ was
 studied by Regensburger, Rosenkranz and Middeke
 \cite{Regensburger-Rosenkranz-Middeke}.


\section{Defining relations for the algebra $\mI_n$}\label{DRAI}

In this section defining relations are found for the algebra
$\mI_n$ and it is proved that the algebra $\mI_n$ is a generalized
Weyl algebra (Proposition \ref{a5Oct9}) of Gelfand-Kirillov
dimension $2n$ (Theorem \ref{5Oct9}) which is neither left nor
right Noetherian (Lemma \ref{a10Oct9}).

{\bf Generalized Weyl Algebras}. Let $D$ be a ring, $\s
=(\s_1,...,\s_n)$ be
 an $n$-tuple  of  commuting ring endomorphisms of $D$,  and
$a=(a_1,...,a_n)$ be  an $n$-tuple  of  elements
 of $D$. 
  The {\it generalized Weyl algebra}
$A=D(\s,a)$ (briefly,  GWA) of degree $n$ is  a ring generated by
$D$  and    $2n$ elements $x_1,...,x_n,$ $y_1,...,$ $y_n$ subject
to the defining relations \cite{Bav-FA91}, \cite{Bav-AlgAnaliz92}:
\begin{align*}
y_ix_i&=a_i ,& x_iy_i&=\s_i(a_i), \\
 x_i d&=\s_i(d)x_i, & d y_i&=y_i\s_i(d), \;\;\;  d \in D,
\end{align*}
$$[x_i,x_j]=[y_i,y_j]=[x_i,y_j]=0, \;\;\;  i\neq j,$$
where $[x, y]=xy-yx$.
  We say that  $a$  and $\s $ are the sets
of {\it defining } elements and endomorphisms of $A$ respectively.
 For a
vector $k=(k_1,...,k_n)\in \mathbb{Z}^n$, let
$v_k=v_{k_1}(1)\cdots v_{k_n}(n)$ where, for $1\leq i\leq n$ and
$m\geq 0$: $\,v_{m}(i)=x_i^m$, $\,v_{-m}(i)=y_i^m$, $v_0(i)=1$. It
follows from  the  definition that $A=\bigoplus_{k\in
\mathbb{Z}^n} A_k$ is a $\mathbb{Z}^n$-graded algebra
($A_kA_e\subseteq A_{k+e},$ for all  $k,e \in \mathbb{Z}^n$),
 where $A_k=v_{k,-}Dv_{k, +}$; $v_{k,+}:=\prod_{k_i>0}v_{k_i}(i)$ and $v_{k,-}=\prod_{k_i<0}v_{k_i}(i)$.
 The tensor product (over the ground  field) $A\t A'$ of generalized Weyl algebras
 of degree $n$ and $n'$ respectively is a GWA of degree $n+n'$:
$$A\otimes A'=D\otimes D'((\tau ,\tau'), (a, a')).$$

Let $\CP_n$ be a polynomial algebra  $K[H_1, \ldots , H_n]$ in $n$
indeterminates  and let $\s =(\s_1,...,\s_n)$ be the  $n$-tuple of
commuting automorphisms of $\CP_n$ such that
 $\s_i(H_i)=H_i-1$ and $\s_i(H_j)=H_j$, for $i\neq j$.  The algebra homomorphism
\begin{equation}\label{AnGWA}
A_n\ra \CP_n ((\s_1,...,\s_n), (H_1, \ldots , H_n)), \;\;
x_i\mapsto  x_i, \; \; \der_i\mapsto y_i, \; \; i=1, \ldots , n,
\end{equation}
  is an isomorphism. We identify the Weyl algebra $A_n$
with the GWA above via this isomorphism. Note that $H_i= \der_i
x_i = x_i\der_i+1$.

It is an experimental fact that many small quantum algebras/groups
are GWAs. More about GWAs and their generalizations the interested
reader can find in \cite{Alev-Far-La-Sol-00, Artam-Cohn-99,
Bav-inform, Bav-Jor-01, Bav-Len-JA-2001-KdimGWA,
Bav-Len-JA-2001-TKMin, Bav-Oyst-1998KdimGWA, Bav-Oyst-2000Adv,
Bav-Oyst-2001TrAMS, Bek-Ben-Fut-04, Cassidy-Shelton,
Far-Sol-SuAl-03, Jordan-Down-up, Hart-06, Kirk-Kuz,
Kir-Mus-Pas-99, Kir-Kuz-05,  Maz-Tur-02, Nauwel-Voyst,
Prest-Pu-02, Rich-Sol-06, Rueda-02,  Staf-sl2-82}.

$\noindent $

Suppose that $A$ is a $K$-algebra that admits two elements $x$ and
$y$ with $yx=1$. The element $xy\in A$ is an idempotent, $(xy)^2=
xy$, and so the set $xyA xy$ is  a $K$-algebra where $xy$ is its
identity element. Consider the linear maps $\s = \s_{x,y}, \tau =
\tau_{x,y}:A \ra A$ which are defined as follows 
\begin{equation}\label{sxyt}
\s (a) =  xay, \;\; \tau (a) = yax.
\end{equation}
Then $\tau \s = \id_A$ and $\s \tau  (a) = xy \cdot a \cdot xy$,
and so the map $\s$ is an algebra monomorphism with $\s (1) = xy$
and
\begin{equation}\label{Asrt}
A= \s (A) \bigoplus \ker (\tau ).
\end{equation}

In more details, $\s (A) \cap \ker (\tau ) =0$ since $\tau \s
=\id_A$. Since $(\s \tau )^2= \s \tau$, we have the equality $A=
\im (\s \tau ) \bigoplus \im (1-\s \tau )$. Clearly, $\im (\s \tau
) \subseteq \im (\s )$ and $\im (1-\s \tau ) \subseteq \ker ( \tau
)$ as $\tau \s = \id_A$. Then $A = \im (\s ) +\ker (\tau )$, i.e.
 (\ref{Asrt}) holds. In general, the map $\tau$ is not an algebra
endomorphism and its kernel is not an ideal of the algebra $A$.
Suppose that the algebra $A$ contains a subalgebra $D$ such that
$\s (D) \subseteq D$ (and so $xy = \s (1) \in D$), and that the
algebra $A$ is generated by $D$, $x$, and $y$. Since $yx=1$, we
have $x^iD_D\simeq D$ and ${}_DDy^i\simeq D$. It follows from the
relations:
\begin{align*}
yx&=1 ,& xy&=\s (1), \\
xd  &=\s (d) x, & d y&=y \s (d), \;\;\;  d \in D,
\end{align*}
that $A= \sum_{i\geq 1} y^iD+\sum_{i\geq 0}Dx^i$. Suppose, in
addition, that the sum is a {\em direct} one. Then the algebra $A$
is the GWA $D(\s , 1)$.

\begin{lemma}\label{a13Apr9}
{\rm \cite{jacaut}}  Keep the assumptions as above, i.e.
$A=D\langle x,y\rangle = \bigoplus_{i\geq 1}
y^iD\bigoplus\bigoplus_{i\geq 0}Dx^i$ and $\s (D)\subseteq D$.
Then $A= D(\s , 1)$. If, in addition, $\tau (D)\subseteq D$ and
the element $xy$ is central in $D$. Then $Dx^i = x^iD$ and $Dy^i =
y^iD$ for all $i\geq 1$.
\end{lemma}

{\it Definition}, \cite{shrekalg}. The 
{\em algebra} $\mathbb{S}_n$ {\em of one-sided inverses} of $P_n$
is an algebra generated over a field $K$  by $2n$ elements $x_1,
\ldots , x_n, y_n, \ldots , y_n$ that satisfy the defining
relations:
$$ y_1x_1=\cdots = y_nx_n=1 , \;\; [x_i, y_j]=[x_i, x_j]= [y_i,y_j]=0
\;\; {\rm for\; all}\; i\neq j,$$ where $[a,b]:= ab-ba$ is  the
algebra  commutator of elements $a$ and $b$.

$\noindent $

By the very definition, the algebra $\mS_n\simeq \mS_1^{\t n}$ is
obtained from the polynomial algebra $P_n$ by adding commuting,
left (but not two-sided) inverses of its canonical generators. The
algebra $\mS_1$ is a well-known primitive algebra
\cite{Jacobson-StrRing}, p. 35, Example 2. Over the field
 $\mathbb{C}$ of complex numbers, the completion of the algebra
 $\mS_1$ is the {\em Toeplitz algebra} which is the
 $\mathbb{C}^*$-algebra generated by a unilateral shift on the
 Hilbert space $l^2(\N )$ (note that $y_1=x_1^*$). The Toeplitz
 algebra is the universal $\mathbb{C}^*$-algebra generated by a
 proper isometry.

The Jacobian algebra $\mA_n$ contains the algebra $\mS_n$ where
$$y_1:=H_1^{-1}\der_1, \ldots, y_n:= H_n^{-1}\der_n.$$
Moreover, the algebra $\mA_n$ is the subalgebra of
 $\End_K(P_n)$ generated by the algebra $\mS_n$ and  the
$2n$ invertible elements $H_1^{\pm 1}, \ldots , H_n^{\pm 1}$ of
$\End_K(P_n)$.

The algebras $\mS_n$ and $\mA_n$ are much more better understood
than the algebra $\mI_n$. We will see that the three classes of
algebras have much in common. In particular, they are GWAs.
Moreover, we will deduce many results for the algebra $\mI_n$ from
known results for the algebras $\mS_n$ and $\mA_n$ in
\cite{Bav-Jacalg}, \cite{shrekalg}, \cite{jacaut}.


{\bf The algebra $\mS_n$ is a GWA}.   Clearly,
$\mathbb{S}_n=\mS_1(1)\t \cdots \t\mS_1(n)\simeq \mathbb{S}_1^{\t
n}$ where $\mS_1(i):=K\langle x_i,y_i \, | \, y_ix_i=1\rangle
\simeq \mS_1$ and $\mS_n=\bigoplus_{\alpha , \beta \in \N^n}
Kx^\alpha y^\beta$ where $x^\alpha := x_1^{\alpha_1} \cdots
x_n^{\alpha_n}$, $\alpha = (\alpha_1, \ldots , \alpha_n)$,
$y^\beta := y_1^{\beta_1} \cdots y_n^{\beta_n}$, $\beta =
(\beta_1,\ldots , \beta_n)$. In particular, the algebra $\mS_n$
contains two polynomial subalgebras $P_n$ and $Q_n:=K[y_1, \ldots
, y_n]$ and is equal, as a vector space,  to their tensor product
$P_n\t Q_n$. Note that also the Weyl algebra $A_n$ is a tensor
product (as a vector space) $P_n\t K[\der_1, \ldots , \der_n]$ of
its two polynomial subalgebras.

When $n=1$, we usually drop the subscript `1' if this does not
lead to confusion (we do the same also for the algebras $A_1$,
$\mA_1$ and $\mI_1$). So, $\mS_1= K\langle x,y\, | \,
yx=1\rangle=\bigoplus_{i,j\geq 0}Kx^iy^j$. For each natural number
$d\geq 1$, let $M_d(K):=\bigoplus_{i,j=0}^{d-1}KE_{ij}$ be the
algebra of $d$-dimensional matrices where $\{ E_{ij}\}$ are the
matrix units, and
$$M_\infty (K) :=
\varinjlim M_d(K)=\bigoplus_{i,j\in \N}KE_{ij}$$ be the algebra
(without 1) of infinite dimensional matrices. The algebra  $
M_\infty (K)=\bigoplus_{k\in \Z}M_\infty (K)_k$ is $\Z$-graded
 where $M_\infty (K)_k:=\bigoplus_{i-j=k}KE_{ij}$ ($M_\infty (K)_kM_\infty
 (K)_l\subseteq M_\infty (K)_{k+l}$ for all $k,l\in \Z$). The algebra $\mS_1$ contains
 the ideal $F:=\bigoplus_{i,j\in
\N}KE_{ij}$, where 
\begin{equation}\label{Eijc}
E_{ij}:= x^iy^j-x^{i+1}y^{j+1}, \;\; i,j\geq 0.
\end{equation}
Note that $E_{ij}=x^iE_{00}y^j$ and $E_{00}=1-xy$.  For all
natural numbers $i$, $j$, $k$, and $l$,
$E_{ij}E_{kl}=\d_{jk}E_{il}$ where $\d_{jk}$ is the Kronecker
delta function.  The ideal $F$ is an algebra (without 1)
isomorphic to the algebra $M_\infty (K)$ via $E_{ij}\mapsto
E_{ij}$. In particular, the algebra $F=\bigoplus_{k\in \Z}F_{1,k}$
is $\Z$-graded where $F_{1,k}:=\bigoplus_{i-j=k}KE_{ij}$
($F_{1,k}F_{1,l}\subseteq F_{1,k+l}$ for all $k,l\in \Z $). For
all $i,j\geq 0$, 
\begin{equation}\label{xyEij}
xE_{ij}=E_{i+1, j}, \;\; yE_{ij} = E_{i-1, j}, \;\; E_{ij}x=E_{i,
j-1}, \;\; E_{ij}y = E_{i, j+1}, \;\;
\end{equation}
where $E_{-1,j}:=0$ and $E_{i,-1}:=0$. 
\begin{equation}\label{xyEij2}
xE_{i, j}=E_{i+1, j+1}x, \;\; E_{ij}y = yE_{i+1, j+1}.
\end{equation}
\begin{equation}\label{mS1d}
\mS_1= K\oplus xK[x]\oplus yK[y]\oplus F,
\end{equation}
the direct sum of vector spaces. Then 
\begin{equation}\label{mS1d1}
\mS_1/F\simeq K[x,x^{-1}]=:L_1, \;\; x\mapsto x, \;\; y \mapsto
x^{-1},
\end{equation}
since $yx=1$, $xy=1-E_{00}$ and $E_{00}\in F$.

$\noindent $

The algebra $\mS_n = \bigotimes_{i=1}^n \mS_1(i)$ contains the
ideal
$$F_n:= F^{\t n }=\bigoplus_{\alpha , \beta \in
\N^n}KE_{\alpha \beta}, \;\; {\rm where}\;\; E_{\alpha
\beta}:=\prod_{i=1}^n E_{\alpha_i \beta_i}(i), \; E_{\alpha_i
\beta_i}(i):= x_i^{\alpha_i}y_i^{\beta_i}-
x_i^{\alpha_i+1}y_i^{\beta_i+1}.$$ Note that $E_{\alpha
\beta}E_{\g \rho}=\d_{\beta \g }E_{\alpha  \rho}$ for all elements
$\alpha, \beta , \g , \rho \in \N^n$ where $\d_{\beta
 \g }$ is the Kronecker delta function.

Using Lemma \ref{a13Apr9}, we can show that the algebra $\mS_1(i)$
is the GWA $\mF_{1,0}(i)(\s_i, 1)$ where
$\mF_{1,0}(i):=K\bigoplus\bigoplus_{k\geq 0}K E_{kk}(i)$ and
$\s_i(a)=x_iay_i$ (moreover, $\s_i(1)=1-E_{00}(i)$ and $ \s_i
(E_{kk}(i))=E_{k+1,k+1}(i)$). Therefore, $\mS_n=
\bigotimes_{i=1}^n \mF_{1,0}(i)(\s_i, 1)=\mF_{n,0}((\s_1, \ldots ,
\s_n), (1, \ldots , 1))$ is a GWA (Lemma 3.3, \cite{jacaut}) where
$\mF_{n,0}:=\bigotimes_{i=1}^n \mF_{1,0}(i)$.
 The algebra $\mF_{n,0}$ is a commutative, non-finitely generated,
 non-Noetherian algebra, it contains the direct sum
 $\bigoplus_{\alpha \in \N^n}KE_{\alpha\alpha}$ of ideals, hence
 $\mF_{n,0}$ is not a prime algebra. The algebra $\mS_n = \bigoplus_{\alpha \in \Z^n} \mS_{n, \alpha}$
is a $\Z^n$-graded algebra where $\mS_{n,
\alpha}=\mF_{n,0}v_\alpha = v_\alpha \mF_{n,0}$ for all $\alpha\in
\Z^n$ where $v_\alpha := \prod_{i=1}^n v_{\alpha_i}(i)$ and
$v_j(i) :=\begin{cases}
x_i^j& \text{if }j\geq 0,\\
y_i^{-j}& \text{if }j<0.\\
\end{cases}$

The map  $\tau_i:\mS_n\ra \mS_n$, $a\mapsto y_iax_i$,  is not an
algebra endomorphism but its restriction to the subalgebra
$\mF_{n,0}$ of $\mS_n$ {\em is} a $K$-algebra epimorphism,
$\tau_i(\mF_{n,0})=\mF_{n,0}$, with $\ker (\tau_i
|_{\mF_{n,0}})=KE_{00}(i)\bigotimes \bigotimes_{j\neq
i}\mF_{1,0}(j)$. For all $j\in \N$ and $d\in \mF_{n,0}$,
$dx_i^j=x_i^j\tau_i^j(d)$ and $ y_i^jd=\tau_i^j(d)y_i^j$.

 The {\em algebra} $\CI_n := K\langle \der_1, \ldots , \der_n , \int_1,
 \ldots , \int_n\rangle$ {\em of integro-differential operators with
 constant coefficients} is canonically isomorphic to the algebra
 $\mS_n$:
\begin{equation}\label{SnIniso}
\mS_n\ra \CI_n, \;\; x_i\mapsto \int_i, \;\; y_i\mapsto
 \der_i, \;\; i=1, \ldots , n.
\end{equation}
For $n=1$ this is obvious since the map above is a well-defined
epimorphism (since $\der \int =1$) which must be an isomorphism as
the algebra $\CI_1$ is non-commutative but any proper factor
algebra of $\mS_1$ is commutative \cite{shrekalg}. Then the
general case follows since $\mS_n \simeq \mS_1^{\t n}$ and $\CI_n
\simeq \CI_1^{\t n}$.


{\bf The Jacobian algebra $\mA_n$ is a GWA}. The Jacobian algebra
$\mA_n = \mD_n ((\s_1, \ldots , \s_n), (1, \ldots , 1))$ is  a GWA
(Lemma 3.4, \cite{jacaut}) where
$\mD_n:=\bigotimes_{i=1}^n\mD_1(i)$,
\begin{eqnarray*}
 \mD_1(i)&=&\CL_1^-(i)\bigoplus \CL_1^+(i)\bigoplus F_{1,0}(i),\;\; F_{1,0}(i):=\bigoplus_{s\geq
0}KE_{ss}(i), \\
\CL_1^-(i)&:=&\bigoplus_{s,t\geq 1}K\frac{1}{(H_i-s)_1^t},\;\;  (H_i-s)_1:= H_i-s+E_{s-1,s-1}(i),\\
 \CL_1^+(i)&:=&K[H^{\pm 1}_i, (H_i+ 1)^{-1},  (H_i+ 2)^{-1}, \ldots ], \;
\end{eqnarray*}
 and $\s_i(a)=x_iay_i$. In more detail, for all
natural numbers $s\geq 0$ and $t\geq 1$,
\begin{eqnarray*}
 \s_i (E_{ss}(i))&=& E_{s+1, s+1}(i), \;\;\;  \s_i ((H_i-t)_1)= (H_i-t-1)_1\s_i (1), \\
 \s_i (H_i)&=&H_i-1=(H_i-1)\s_i (1)
= (H_i-1)_1\s_i (1).
\end{eqnarray*}
 The algebra $\mD_n$ is a commutative, non-finitely generated,
 non-Noetherian algebra, it contains the direct sum
 $\bigoplus_{\alpha \in \N^n}KE_{\alpha\alpha}$ of ideals (and so
 $\mD_n$ is not a prime algebra). Note that $H_i E_{ss}(i)= E_{ss}(i) H_i = (s+1)E_{ss}(i)$.
Clearly, $\mA_n = \bigotimes_{i=1}^n
\mA_1(i)=\bigotimes_{i=1}^n\mD_1(i)(\s_i, 1)$.  The algebra $\mA_n
= \bigoplus_{\alpha \in \Z^n} \mA_{n, \alpha}$ is a $\Z^n$-graded
algebra where $\mA_{n,\alpha}=\mD_n v_\alpha = v_\alpha \mD_n$,
$v_\alpha := \prod_{i=1}^n v_{\alpha_i}(i)$ and 
\begin{equation}\label{Anvji}
v_j(i) :=
\begin{cases}
x_i^j& \text{if }j\geq 0,\\
y_i^{-j}& \text{if }j<0.\\
\end{cases}
\end{equation}

The map $\tau_i:\mA_n\ra \mA_n$, $a\mapsto y_iax_i$, is not an
algebra endomorphism but its restriction to the subalgebra $\mD_n$
of $\mA_n$ {\em is} a $K$-algebra epimorphism,
$\tau_i(\mD_n)=\mD_n$,  $\ker (\tau_i
|_{\mD_n})=KE_{00}(i)\bigotimes\bigotimes_{j\neq i}\mD_1(j)$. In
more detail, for $a,b\in \mD_1$,
$$\tau (a) \tau (b) = ya(1-E_{00})bx= \tau (ab)-yaE_{00}bx= \tau
(ab ),$$ since $aE_{00}b \in KE_{00}$ and $yE_{00}=0$. For all
$j\in \N$ and $d\in \mD_n$, $dx_i^j=x_i^j\tau_i^j(d)$ and $
y_i^jd=\tau_i^j(d)y_i^j$. Indeed, when $n=1$, $x\tau (d) =
(1-E_{00})dx=dx-E_{00}dx=dx$ since $E_{00}d\in KE_{00}$ and
$E_{00}x=0$.

{\bf The algebra $\mI_n$ is a GWA}.  Since $x_i= \int_iH_i$, the
algebra $\mI_n$ is generated by the elements $\{ \der_i , H_i,
\int_i\, | \, i=1, \ldots , n\}$, and $\mI_n =\bigoplus_{i=1}^n
\mI_1(i)$ where $\mI_1(i):= K\langle \der_i , H_i,
\int_i\rangle=K\langle \der_i , x_i, \int_i\rangle\simeq \mI_1$.
By (\ref{SnIniso}), when $n=1$ the following elements of the
algebra $\mI_1=K\langle \der , H, \int \rangle$,
\begin{equation}\label{eijdef}
e_{ij}:=\int^i\der^j-\int^{i+1}\der^{j+1}, \;\; i,j\in \N ,
\end{equation}
satisfy the relations: $e_{ij}e_{kl}=\d_{jk}e_{il}$. Note that
$e_{ij}=\int^ie_{00}\der^j$. The matrices of the linear maps
$e_{ij}\in \End_K(K[x])$ with respect to the basis $\{ x^{[s]}:=
\frac{x^s}{s!}\}_{s\in \N}$ of the polynomial algebra $K[x]$  are
the elementary matrices, i.e.
$$ e_{ij}*x^{[s]}=\begin{cases}
x^{[i]}& \text{if }j=s,\\
0& \text{if }j\neq s.\\
\end{cases}$$
It follows that 
\begin{equation}\label{eijEij}
e_{ij}=\frac{j!}{i!}E_{ij},
\end{equation}
 $Ke_{ij}=KE_{ij}$, and
$F=\bigoplus_{i,j\geq 0}Ke_{ij}\simeq M_\infty (K)$. Moreover,
$F_n=\bigoplus_{\alpha, \beta \in \N^n} Ke_{\alpha\beta}$ where
$e_{\alpha\beta}:= \prod_{i=1}^n e_{\alpha_i\beta_i}(i)$ and
$e_{\alpha_i\beta_i}(i):=\int_i^{\alpha_i}\der_i^{\beta_i}-\int_i^{\alpha_i+1}\der_i^{\beta_i+1}$.

 The next proposition gives a finite set of defining relations for
 the algebra $\mI_n$ and shows that the algebra $\mI_n$ is a GWA
 (and so we have another set of defining relations for the
 algebra $\mI_n$).
\begin{proposition}\label{a5Oct9}
\begin{enumerate}
\item The algebra $\mI_n$ is generated by the elements $\{ \der_i,
\int_i, H_i\, | \, i=1, \ldots , n\}$ that satisfy the  defining
relations:
\begin{eqnarray*}
\forall i& : & \der_i\int_i= 1, \;\; [H_i, \int_i]=\int_i, \;\;
[ H_i, \der_i]=-\der_i, \;\; H_i (1-\int_i\der_i) = (1-\int_i\der_i) H_i=1-\int_i\der_i, \\
\forall i\neq j &:   & a_ia_j= a_ja_i \;\; {\rm where}\;\; a_k\in \{ \der_k, \int_k, H_k\}. \\
\end{eqnarray*}
\item The algebra $\mI_n= \bigotimes_{i=1}^n D_1(i)(\s_i,
1)=D_n((\s_1, \ldots , \s_n), (1,\ldots , 1))$ is a GWA
($\int_i\lra x_i$, $\der_i\lra y_i$, $H_i\lra H_i$) where
$D_n:=\bigotimes_{i=1}^n D_1(i)$, $D_1(i):=K[H_i]\bigoplus
\bigoplus_{j\geq 0}Ke_{jj}(i)$, $H_ie_{jj}(i) = e_{jj}(i)H_i=
(j+1)e_{jj}(i)$, and the $K$-algebra endomorphisms $\s_i$  of
$D_n$ are given by the rule $\s_i (a):= \int_ia\der_i$ ($\s_i(H_i)
= H_i-1$, $\s_ie_{jj}(i))=e_{j+1, j+1}(i)$). Moreover, the algebra
$\mI_n=\bigoplus_{\alpha \in \Z^n}\mI_{n, \alpha}$ is
$\Z^n$-graded where $\mI_{n,\alpha}= D_n v_\alpha = v_\alpha D_n $
for all $\alpha \in \Z^n$ where $v_\alpha:=
\prod_{i=1}^nv_{\alpha_i}(i)$ and $v_j(i):=\begin{cases}
\int_i^j& \text{if } j>0,\\
1& \text{if } j=0.\\
\der_i^{-j} & \text{if } j<0.\\
\end{cases} $
\item {\rm (The canonical basis for the algebra $\mI_n$)} $\mI_n =
\bigoplus_{\alpha \in \Z^n}\mI_{n,\alpha}$ and, for all $\alpha
\in \Z^n$, $\mI_{n,\alpha}=v_{\alpha , +}D_n v_{\alpha , -}\simeq
D_n$ ($v_{\alpha , +}d v_{\alpha , -}\lra d$) where $v_{\alpha ,
+}:=\prod_{\alpha_i>0} v_{\alpha_i}(i)$ and $v_{\alpha,
-}:=\prod_{\alpha_i<0}v_{\alpha_i}(i)$. So, each element $a\in
\mI_n$ is a unique finite sum $a=\sum_{\alpha \in \Z^n}v_{\alpha ,
+}a_\alpha  v_{\alpha , -}$ for unique elements $a_\alpha \in
D_n$.
\end{enumerate}
\end{proposition}

{\it Proof}. It suffices to prove the statements for $n=1$ since
$\mI_n= \bigotimes_{i=1}^n \mI_1(i)$. So, let $n=1$ and $\mI_1'$
be an algebra generated by symbols $\der$, $\int$, and $H$ that
satisfy the defining relations of statement 1. The algebra $\mI_1$
is generated by the elements  $\der$, $\int$, and $H$; and they
satisfy the defining relations of statement 1 as we can easily
verify. Therefore, there is the natural algebra epimorphism
$\mI_1'\ra \mI_1$ given by the rule: $\der \mapsto \der$,
$\int\mapsto \int$, $H\mapsto H$. It follows from the relations of
statement 1 and from the equalities $e_{ij}=\begin{cases}
\int^{i-j}e_{jj}& \text{if }i\geq j,\\
e_{ii}\der^{j-i}& \text{if } i<j, \\
\end{cases}$ that
$$ \mI_1'=\sum_{i\geq 1}D_1'\der^i+D_1'+\sum_{i\geq 1}\int^iD_1'
=\sum_{i\geq 1}\der^iD_1'+D_1'+\sum_{i\geq 1}D_1'\int^i$$ where
  $D_1':=K\langle H\rangle +\sum_{i\geq
0}Ke_{ii}$. Since $\der \int =1$,  the left $D_1'$-modules  $D_1'$
and $ D_1'\der^i$ are isomorphic, and the right $D_1'$-modules
$\int^iD_1'$ and $D_1'$ are isomorphic.
 Using the $\Z$-grading of the Jacobian
algebra $\mA_1$ and the fact that $\mI_1\subseteq \mA_1$, we have
$$\mI_1=\bigoplus_{i\geq 1}D_1\der^i\bigoplus D_1\bigoplus
\bigoplus_{i\geq 1} \int^iD_1=\bigoplus_{i\geq
1}\der^iD_1\bigoplus D_1\bigoplus \bigoplus_{i\geq 1} D_1\int^i$$
where $D_1= K[H]\bigoplus \bigoplus_{i\geq 0}Ke_{ii}=
K[H]\bigoplus \bigoplus_{i\geq 0}KE_{ii}$ since $Ke_{ii}=KE_{ii}$.
Note that the left $D_1$-modules $D_1$ and $D_1\der^i$ are
isomorphic and the right $D_1$-modules $D_1$ and $\int^iD_1$ are
isomorphic since $\der \int =1$.  This implies that the sum for
$\mI_1'$ above is a direct one. Therefore, $\mI_1'\simeq \mI_1$
and the relations in statement 1 are defining relations for the
algebra $\mI_1$ and $D_1'=D_1$. The condition of Lemma
\ref{a13Apr9} hold, and so $\mI_1=D_1(\s , 1)$ with
$D_1\int^i=\int^iD_1$ and $D_1\der^i= \der^i D_1$ for all $i\geq
1$. The proof of  statements 1 and 2 of the  proposition is
complete.

Statement 3 follows from statement 2 and the fact that, for all
$\alpha \in \Z^n$, the linear map $\mI_{n, \alpha}\ra D_n$,
$b\mapsto u_{\alpha, -}bu_{\alpha, +}$, is a bijection since
$u_{\alpha , -}v_{\alpha , +}=1$ and $v_{\alpha, -}u_{\alpha ,
+}=1$ where $u_{\alpha, -}:= \prod_{\alpha_i>0}v_{-\alpha_i}(i)$
and $ u_{\alpha, +}:= \prod_{\alpha_i<0}v_{-\alpha_i}(i)$.  $\Box
$

$\noindent $

{\it Definition}. For each element $a\in \mI_n$, the unique sum
for $a$ in statement 3 of Proposition \ref{a5Oct9} is called the
{\em canonical form} of $a$.

$\noindent $

The map $\tau_i :\mI_n\ra \mI_n$, $a\mapsto \der_i a\int_i$, is
not an algebra endomorphism but its restriction to the subalgebra
$D_n$ of  $\mI_n$ is a $K$-algebra {\em epimorphism}, $\tau_i(D_n)
= D_n$ with $\ker (\tau_i|_{D_n}) =
Ke_{00}(i)\bigoplus\bigoplus_{j\neq i}D_1(j)$. In more detail, for
$n=1$ and $a,b\in D_1$, $$\tau (a) \tau (b) = \der
a(1-e_{00})b\int= \tau (ab) - \der a e_{00}b\int= \tau (ab)$$
since $ ae_{00}b\in Ke_{00}$ and $\der e_{00}=0$. For all $j\in
\N$ and $d\in D_n$, $d \int_i^j=\int_i^j\tau_i^j(d)$ and $\der_i^j
d= \tau_i^j(d) \der_i^j$. Indeed, for $n=1$, $\int \tau (d) =
(1-e_{00})d\int= d\int - e_{00}d\int= d\int$ since $e_{00}d\in
Ke_{00}$ and $e_{00}\int =0$. Note that 
\begin{equation}\label{tiHj1}
\tau_i(H_j) = H_j+\d_{ij}\;\; {\rm and}\;\;
\tau_i(e_{st}(j))=\begin{cases}
e_{s-1, t-1}(i)& \text{if }i=j,\\
e_{st}(j)& \text{if }i\neq j.\\
\end{cases}
\end{equation}
It follows that
 $\bigcap_{i=1}^n\ker(\tau_i|_{D_n}) =
K\prod_{i=1}^ne_{00}(i)= Ke_{00}$ and  $\bigcap_{i=1}^n
\ker(\tau_i|_{D_n}-1) = K$.

For the definition and properties of the Gelfand-Kirillov
dimension $\GK $ the reader is referred to \cite{KL} and
\cite{MR}.

\begin{theorem}\label{5Oct9}
The Gelfand-Kirillov dimension $\GK (\mI_n)$ of the algebra
$\mI_n$ is $2n$.
\end{theorem}

{\it Proof}. Since $A_n\subseteq \mI_n$, we have the inequality
$2n=\GK (A_n) \leq \GK (\mI_n)$. To prove the reverse inequality
let  us consider the standard filtration $\{ \mI_{n,i}\}_{i\in
\N}$ of the algebra $\mI_n$ with respect to the set of generators
$\{ \der_i , H_i, \int_i\, | \, i=1, \ldots , n\}$ of the algebra
$\mI_n$. By Proposition \ref{5Oct9}, $\mI_{n,i}\subseteq
\mI_{n,i}':= \bigoplus_{|\alpha |\leq i} v_\alpha D_{n,i}$ where
$D_{n,i}:= \bigotimes_{j=1}^n D_{1,i}(j)$ and $D_{1,i}(j) :=
\bigoplus_{s=0}^iKH_j^s\bigoplus\bigoplus_{t=0}^iKe_{tt}(j)$. Then
$\dim (\mI_{n,i})\leq \dim (\mI_{n,i}')\leq (2i+1)^n (2i+2)^n$,
and so $\GK(\mI_n)\leq 2n$, as required. $\Box $

\begin{lemma}\label{a10Oct9}
The algebra $\mI_n$ is neither left nor right Noetherian.
 Moreover, it contains infinite direct sums of nonzero left (resp.
right) ideals.
\end{lemma}

{\it Proof}. Since $\mI_n \simeq \mI_1^{\t n}$, it suffices to
prove the lemma for $n=1$. The ideal $F=\bigoplus_{i,j\geq
0}KE_{ij}$ of the algebra $\mI_1$ is the infinite direct sum
 $\bigoplus_{j\geq 0}(\bigoplus_{i\geq 0}KE_{ij})$ (resp.
$\bigoplus_{i\geq 0}(\bigoplus_{j\geq 0}KE_{ij})$) of nonzero left
(resp. right) ideals, and the statements follow. $\Box $


\section{Ideals of  the algebra $\mI_n$}\label{IOTAI}

In this section, we prove that the restriction map (Theorem
\ref{7Oct9}) from the set of ideals $\CJ (\mA_n)$ of the algebra
$\mA_n$ to the set of  ideals $\CJ (\mI_n)$ of the algebra $\mI_n$
is a bijection that respects the three operations on ideals: sum,
intersection and product. As a consequence, we obtain many results
 for the ideals of the algebra $\mI_n$ using similar results for the
ideals of the algebra $\mA_n$ in  \cite{Bav-Jacalg}, see Corollary
\ref{b10Oct9} and Corollary \ref{c10Oct9}: a classification of all
the ideals of $\mI_n$ (there are only finitely many of them) and a
classification of prime ideals of $\mI_n$, etc.

$\noindent $

{\it Definition}. Let $A$ and $B$ be algebras, and let $\CJ (A)$
and $\CJ (B)$ be their lattices of ideals. We say that a bijection
$f: \CJ (A) \ra \CJ (B)$ is an {\em isomorphism} if $f(\ga *\gb )
= f(\ga )*f(\gb )$  for $*\in \{ +, \cdot, \cap \}$, and in this
case we say that the algebras $A$ and $B$ are {\em ideal
equivalent}. The ideal equivalence is an equivalence relation on
the class of algebras.

The next theorem shows that the algebras $\mA_n$ and $\mI_n$ are
ideal equivalent.
\begin{theorem}\label{7Oct9}
The restriction map $\CJ ( \mA_n) \ra \CJ (\mI_n)$, $\ga \mapsto
\ga^r:= \ga \cap \mI_n$, is an isomorphism (i.e. $(\ga_1*\ga_2)^r=
\ga_1^r*\ga_2^r$ for $*\in \{ +, \cdot, \cap \}$) and its inverse
is the extension map $\gb\mapsto \gb^e:= \mA_n \gb \mA_n$.
\end{theorem}

{\it Proof}. The theorem follows from Theorem \ref{10Oct9}. $\Box
$

$\noindent $

Recall that $\mF_{n, 0}\subset \mI_n\subset \mA_n\subset
\End_K(P_n)$. The subset of $\CJ (\mF_{n, 0})$, $\CJ (\mF_{n,
0})_{\s, \tau}:=\{ \gb \in \CJ (\mF_{n, 0})\, | \, \s_i (\gb )
\subseteq \gb , \tau_i (\gb ) \subseteq \gb$ for all $i=1, \ldots
, n\}$, is closed under addition, multiplication and intersection
of ideals where $\s_i(a) = \int_ia\der_i$ and $\tau_i (a) =
\der_ia\int_i$ (recall that the maps $\s_i, \tau_i : \mF_{n, 0}\ra
\mF_{n, 0}$ are $K$-algebra homomorphisms; $\tau_i(1)=1$ but
$\s_i(1) = 1-e_{00}(i)$).

\begin{theorem}\label{10Oct9}
\begin{enumerate}
\item The restriction map $\CJ ( \mI_n) \ra \CJ (\mF_{n, 0})_{\s
,\tau} $, $\ga \mapsto \ga^r:= \ga \cap \mF_{n, 0}$, is an
isomorphism (i.e. $(\ga_1*\ga_2)^r= \ga_1^r*\ga_2^r$ for $*\in \{
+, \cdot, \cap \}$) and its inverse is the extension map
$\gb\mapsto \gb^e:= \mI_n \gb \mI_n$. \item The restriction map
$\CJ ( \mA_n) \ra \CJ (\mF_{n, 0})_{\s ,\tau}$, $\ga \mapsto
\ga^r:= \ga \cap \mF_{n, 0}$, is an isomorphism (i.e.
$(\ga_1*\ga_2)^r= \ga_1^r*\ga_2^r$ for $*\in \{ +, \cdot, \cap
\}$) and its inverse is the extension map $\gb\mapsto \gb^e:=
\mA_n \gb \mA_n$.
\end{enumerate}
\end{theorem}

The proof of Theorem \ref{10Oct9} is given at the end of this
section. Now, we obtain some consequences of Theorem \ref{7Oct9}.

 The next corollary shows that the ideal
theory of $\mI_n$ is `very arithmetic.' Let $\CB_n$ be the set of
all functions $f:\{ 1, 2, \ldots , n\} \ra  \{ 0,1\}$. For each
function $f\in \CB_n$, $I_f:= I_{f(1)}\t \cdots \t I_{f(n)}$ is
the ideal of $\mI_n$ where $I_0:=F$ and $I_1:= \mI_1$.  Let
$\CC_n$ be the set of all  subsets of $\CB_n$ all distinct
elements of which are incomparable (two distinct elements $f$ and
$g$ of $\CB_n$ are {\em incomparable} if neither $f(i)\leq  g(i)$
nor $f(i)\geq g(i)$ for all $i$). For each $C\in \CC_n$, let
$I_C:= \sum_{f\in C}I_f$, the ideal of $\mI_n$. The number $\gd_n$
of elements in the set $\CC_n$ is called the {\em Dedekind
number}. It appeared in the paper of Dedekind
\cite{Dedekind-1871}. An asymptotic of the Dedekind numbers was
found by Korshunov \cite{Korshunov-1977}.

\begin{corollary}\label{b10Oct9}
\begin{enumerate}
\item The algebra $\mI_n$ is a prime algebra. \item The set of
height one prime ideals of the algebra $\mI_n$ is $\{ \gp_1:=
F\t\mI_{n-1}, \gp_1:= \mI_1\t F\t\mI_{n-2},\ldots , \gp_n:=
\mI_{n-1}\t F\}$. \item Each ideal of the algebra $\mI_n$ is an
idempotent ideal ($\ga^2= \ga$). \item The ideals of the algebra
$\mI_n$ commute ($\ga \gb = \gb \ga$). \item The lattice $\CJ
(\mI_n)$ of ideals of the algebra $\mI_n$ is distributive. \item
The classical Krull dimension $\clKdim (\mI_n)$ of the algebra
$\mI_n$ is $n$. \item $\ga \gb = \ga \cap \gb$ for all ideals $\ga
$ and $\gb $ of the algebra $\mI_n$. \item The ideal $\ga_n :=
\gp_1+\cdots + \gp_n$ is the largest (hence, the only maximal)
ideal of $\mI_n$ distinct from $\mI_n$, and $F_n = F^{\t
n}=\bigcap_{i=1}^n \gp_i$ is the smallest nonzero ideal of
$\mI_n$. \item {\rm (A classification of  ideals of $\mI_n$)} The
map
 $\CC_n\ra \CJ (\mI_n)$, $C\mapsto I_C:= \sum_{f\in C}I_f$
 is a bijection where $I_\emptyset :=0$. The number of ideals of
 $\mI_n$ is the Dedekind number $\gd_n$.  Moreover,
$2-n+\sum_{i=1}^n2^{n\choose i}\leq \gd_n \leq 2^{2^n}$. For
$n=1$, $F$ is the unique proper ideal of the algebra $\mI_1$.
\item {\rm (A classification of prime ideals of $\mI_n$)} Let
$\Sub_n$ be the set of all subsets of $\{ 1, \ldots , n\}$. The
map $\Sub_n\ra \Spec (\mI_n)$, $ I\mapsto \gp_I:= \sum_{i\in
I}\gp_i$, $\emptyset \mapsto 0$, is a bijection, i.e. any nonzero
prime ideal of $\mI_n$ is a unique sum of primes of height 1;
$|\Spec (\mI_n)|=2^n$; the height of $\gp_I$ is $| I|$; and
 \item  $\gp_I\subset \gp_J$ iff $I\subset
 J$.
\end{enumerate}
\end{corollary}

{\it Proof}. By Theorem \ref{7Oct9}, the statements hold for the
algebra $\mI_n$ since the same statements hold for the algebra
$\mA_n$, and below references are given for their proofs in
\cite{Bav-Jacalg}.

1. Corollary 2.7.(5).

2.  Corollary 3.5.

3.  Theorem 3.1.(2).

4.  Corollary 3.10.(3).

5. Theorem 3.11.

6. Corollary 3.7.

7. Corollary 3.10.(3).

8. Corollary 2.7.(4,7).

9. Theorem 3.1.

10. Corollary 3.5.

11. Corollary 3.6.  $\Box $

$\noindent $

For each ideal $\ga$ of $\mI_n$, $\Min (\ga )$ denotes the set of
minimal primes over $\ga$. Two distinct prime ideals $\gp $ and
$\gq$ are called {\em incomparable} if neither $\gp \subseteq \gq$
nor $\gp\supseteq \gq$. The algebras $\mI_n$ have  beautiful ideal
theory  as the following unique factorization properties
demonstrate.

\begin{corollary}\label{c10Oct9}
\begin{enumerate}
\item Each ideal $\ga$ of $\mI_n$ such that $\ga \neq \mI_n$ is a
unique product of incomparable primes, i.e. if $\ga = \gq_1\cdots
\gq_s= \gr_1\cdots \gr_t$ are two such products then $s=t$ and
$\gq_1= \gr_{\s (1)}, \ldots , \gq_s= \gr_{\s (s)}$ for a
permutation $\s$ of $\{ 1, \ldots, n\}$. \item Each ideal $\ga$ of
$\mI_n$ such that $\ga \neq \mI_n$ is a unique intersection of
incomparable primes, i.e. if $\ga = \gq_1\cap \cdots\cap  \gq_s=
\gr_1\cap \cdots \cap \gr_t$ are two such intersections then $s=t$
and $\gq_1= \gr_{\s (1)}, \ldots , \gq_s= \gr_{\s (s)}$ for a
permutation $\s$ of $\{ 1, \ldots, n\}$. \item For each ideal
$\ga$ of $\mI_n$ such that $\ga \neq \mI_n$, the sets of
incomparable primes in statements 1 and 2 are the same, and so
$\ga=\gq_1\cdots \gq_s=\gq_1\cap \cdots\cap \gq_s$. \item The
ideals $\gq_1,\ldots , \gq_s$ in statement 3 are the minimal
primes of $\ga$, and so  $\ga = \prod_{\gp \in \Min (\ga )}\gp
=\cap_{\gp \in \Min (\ga )}\gp$.
\end{enumerate}
\end{corollary}

{\it Proof}. The same statements are true for the algebra $\mA_n$
(Theorem 3.8, \cite{Bav-Jacalg}).  Now, the corollary follows from
Theorem \ref{7Oct9}. $\Box $

$\noindent $

The next corollary gives all decompositions of an ideal as a
product or intersection of ideals.
\begin{corollary}\label{d10Oct9}

 Let $\ga$ be an ideal of $\mI_n$, and $\CM$ be the minimal
elements with respect to inclusion of the set of minimal primes
of a set of ideals $\ga_1, \ldots , \ga_k$ of $\mI_n$. Then

\begin{enumerate}
\item $\ga = \ga_1\cdots \ga_k$  iff $\Min (\ga ) = \CM$.\item
$\ga = \ga_1\cap \cdots\cap \ga_k$  iff $\Min (\ga ) = \CM$.
\end{enumerate}\end{corollary}

{\it Proof}. The same statements are true for the algebra $\mA_n$
(Theorem 3.12, \cite{Bav-Jacalg}), and the corollary follows from
Theorem \ref{7Oct9}. $\Box $

This is a rare example of a noncommutative  algebra of classical
Krull dimension $>1$ where one has a complete picture of
decompositions of ideals. Recall that a ring $R$ of finite
classical Krull dimension is called {\em catenary} if, for each
pair of prime ideals $\gp$ and $\gq$ with $\gp\subseteq \gq$, all
maximal chains of prime ideals, $\gp_0=\gp\subset \gp_1\subset
\cdots \subset \gp_l=\gq$, have the same length $l$.

\begin{corollary}\label{a12Oct9}
The algebra $\mI_n$ is catenary.
\end{corollary}

{\it Proof}. This follows from Corollary \ref{b10Oct9}.(10, 11).
$\Box $

\begin{corollary}\label{e10Oct9}
 The same statements (with obvious modifications) of Corollaries
 \ref{b10Oct9} and \ref{c10Oct9} hold for the ideals $\CJ
 (\mF_{n,0})_{\s , \tau }$ of the algebra $\mF_{n,0}$ rather than $\mI_n$
  (we leave it to the reader to formulate them).
\end{corollary}

\begin{proposition}\label{B17Oct9}
The polynomial algebra $P_n$ is the only (up to isomorphism)
faithful simple $\mI_n$-module.
\end{proposition}

{\it Proof}. The $\mI_n$-module $P_n$ is faithful (as
$\mI_n\subset \End_K(P_n)$) and simple since the $A_n$-module
$P_n$ is simple and  $A_n\subset \mI_n$. Let $M$ be a faithful
simple $\mI_n$-module. Then $F_nM\neq 0$, i.e. $e_{0\beta}m\neq 0$
for some elements $\beta \in \N^n$ and $m\in M$. The
$\mI_n$-module $P_n\simeq \mI_n/\sum_{i=1}^n \mI_n \der_i$ is
simple. Therefore, the $\mI_n$-module epimorphism $P_n\ra M= \mI_n
e_{0\beta}m = \sum_{\alpha \in \N^n} Ke_{\alpha\beta}m$, $1\mapsto
e_{0\beta}m$, is an isomorphism. The proof of the proposition is
complete.  $\Box $

For a ring $R$ and its element $r\in R$, $\ad (r) : s\mapsto
[r,s]:=rs-sr$ is the {\em inner derivation} of the ring $R$
associated with the element $r$.

{\bf The proof of Theorem \ref{10Oct9}}.  We split the proof of
Theorem \ref{10Oct9} into several statements (which are
interesting on their own) to  make the proof clearer.

The algebra $D_n$ is the zero component of the $\Z^n$-graded
algebra $\mI_n=D_n((\s_1, \ldots , \s_n), (1, \ldots , 1))$, hence
$\s_i(D_n) \subseteq D_n $ and  $\tau_i(D_n) \subseteq D_n$ for
all $i$ where $\s_i(a) = \int_ia\der_i$ and $\tau_i (a) =
\der_ia\int_i$. Let $\CJ (D_n)_{\s , \tau}:= \{ \gb \in \CJ (D_n)
\, | \, \s_i (\gb )\subseteq  \gb , \tau_i (\gb ) \subseteq \gb$
for all $i=1, \ldots , n\}$. Similarly, the algebra $\mD_n$ is the
zero component of the $\Z^n$-graded algebra $\mA_n=\mD_n((\s_1,
\ldots , \s_n), (1, \ldots , 1))$, hence $\s_i(\mD_n) \subseteq
\mD_n $ and $\tau_i(\mD_n) \subseteq \mD_n$ for all $i$ where
$\s_i(a) = x_iay_i$ and $\tau_i (a) = y_iax_i$. Let $\CJ
(\mD_n)_{\s , \tau}:= \{ \gb \in \CJ (\mD_n) \, | \, \s_i (\gb )
\subseteq \gb , \tau_i (\gb ) \subseteq \gb$ for all $i=1, \ldots
, n\}$.

\begin{theorem}\label{A10Oct9}
\begin{enumerate}
\item For each ideal $\ga$ of the algebra $\mI_n$, $\ga =
\bigoplus_{\alpha \in \Z^n}v_{\alpha, +}\ga^r v_{\alpha, -}$ where
$\ga^r:= \ga \cap D_n \in \CJ (D_n)_{\s, \tau }$ and, for each
ideal $\gb \in \CJ (D_n)_{\s , \tau}$, $\gb^e:= \mI_n \gb
\mI_n=\bigoplus_{\alpha \in \Z^n}v_{\alpha, +}\gb  v_{\alpha, -}$
 where $v_{\alpha, +}:=\prod_{\alpha_i>0}v_{\alpha_i}(i)$,
 $v_{\alpha, -}:=\prod_{\alpha_i<0}v_{\alpha_i}(i)$, and $v_j(i)$
 is defined in Proposition \ref{a5Oct9}.(2).
\item For each ideal $\ga$ of the algebra $\mA_n$, $\ga =
\bigoplus_{\alpha \in \Z^n}v_{\alpha, +}\ga^r v_{\alpha, -}$ where
$\ga^r:= \ga \cap \mD_n \in \CJ (\mD_n)_{\s, \tau }$ and, for each
ideal $\gb \in \CJ (\mD_n)_{\s , \tau}$, $\gb^e:= \mA_n \gb
\mA_n=\bigoplus_{\alpha \in \Z^n}v_{\alpha, +}\gb  v_{\alpha, -}$
 where $v_{\alpha, \pm}$ are as above but  the elements $v_j(i)$
 are defined in (\ref{Anvji}).
\end{enumerate}
\end{theorem}

{\it Proof}. 1. Let $\ga$ be an ideal of the algebra $\mI_n$. The
algebra $\mI_n=\bigoplus_{\alpha \in \Z^n}\mI_{n,\alpha}$ is a
$\Z^n$-graded algebra  with $\mI_{n,\alpha}:= \bigcap_{i=1}^n \ker
( \ad (H_i) -\alpha_i)$ for all $\alpha \in \Z^n$. Then $\ga$ is a
homogeneous ideal, that is $\ga = \bigoplus_{\alpha \in \Z^n}
\ga_\alpha$ where $\ga_\alpha := \ga \cap \mI_{n, \alpha}$. The
ideal $\ga_0:= \ga \cap D_n = \ga^r$ of the algebra $D_n$ belongs
to the set $\CJ (D_n)_{\s , \tau }$ since $\s_i ( \ga_0)= \int_i
\ga_0\der_i \subseteq \ga_0$ and $\tau_i ( \ga_0) = \der_i
\ga_0\int_i \subseteq \ga_0$ for all $i=1, \ldots , n$. By
Proposition \ref{a5Oct9}.(2), $\ga_\alpha =v_{\alpha, +}\gb_\alpha
v_{\alpha, -}$ for some ideal $\gb_\alpha$ of the algebra $D_n$:
$\ga_\alpha = D_n \ga_\alpha D_n = D_nv_{\alpha, +}\gb_\alpha
v_{\alpha, -}D_n = v_{\alpha, +}\tau_{\alpha, +}(D_n)\gb_\alpha
\tau_{\alpha, -}(D_n) v_{\alpha, -}= v_{\alpha, +}D_n\gb_\alpha
D_n v_{\alpha, -}$ since $\tau_{\alpha, \pm }(D_n) = D_n$ where
$\tau_{\alpha, +}:= \prod_{\alpha_i>0}\tau_i$ and $\tau_{\alpha,
-}:= \prod_{\alpha_i<0}\tau_i$. Let 
\begin{equation}\label{uua}
u_{\alpha, -}:= \prod_{\alpha_i>0}v_{-\alpha_i}(i), \;\;
u_{\alpha, +}:= \prod_{\alpha_i<0}v_{-\alpha_i}(i).
\end{equation}
 The ideal $\gb_\alpha$ is
unique since $v_{\alpha, +}\gb_\alpha v_{\alpha, -}=v_{\alpha,
+}\gb_\alpha' v_{\alpha, -}$ implies
$$ \gb_\alpha = 1\cdot \gb_\alpha \cdot 1=u_{\alpha, -} v_{\alpha, +}\gb_\alpha
v_{\alpha, -}u_{\alpha, +}=u_{\alpha, -} v_{\alpha, +}\gb_\alpha'
v_{\alpha, -}u_{\alpha, +}=1\cdot \gb_\alpha' \cdot
1=\gb_\alpha'.$$ Moreover, $\gb_\alpha = \ga_0$ for all $\alpha
\in \Z^n$ since $\ga_0\supseteq u_{\alpha, -}\ga_{n, \alpha }
u_{\alpha, +}=u_{\alpha, -} v_{\alpha, +}\gb_\alpha v_{\alpha,
-}u_{\alpha, +}=1\cdot \gb_\alpha \cdot 1 =\gb_\alpha$. On the
other hand, $\ga_{n, \alpha} \supseteq v_{\alpha, +}\ga_0
v_{\alpha, -}$, and so $\ga_0\subseteq \gb_\alpha$.

Let $\gb \in \CJ (D_n)_{\s , \tau}$. Then $\gb^{er}=\gb$ since
\begin{eqnarray*}
 \gb &\subseteq & \gb^{er}=(\mI_n \gb \mI_n)^r=
 \sum_{\alpha \in \Z^n}\mI_{n, \alpha} \gb \mI_{n, -\alpha} \\
 &=& \sum_{\alpha \in \Z^n} v_\alpha D_n \gb D_n v_{-\alpha} \;\;\;
 {\rm (Proposition \;\ref{a5Oct9}.(2))}\\
 &=& \sum_{\alpha \in \Z^n}\s_{\alpha, +}\tau_{\alpha, -}(\gb )
 v_\alpha v_{-\alpha }\subseteq \gb D_n = \gb ,
\end{eqnarray*}
where $\s_{\alpha, +}:=\prod_{\alpha_i>0}\s_i$.  Therefore, $\gb^e
= \bigoplus_{\alpha \in \Z^n} v_{\alpha, +} \gb v_{\alpha, -}$, by
the first part of statement 1.

2. Repeat the proof of statement 1 replacing $(\mI_n , D_n)$ by
 $(\mA_n , \mD_n)$ and making obvious modifications. $\Box $

$\noindent $

For each function $f\in \CB_n$, let

\begin{eqnarray*}
 \gf_f&:=& \gf_{f(1)}\t \cdots \t \gf_{f(n)}\;\;
 {\rm where}\;\; \gf_0:=F_{1,0}, \;\; \gf_1:= \mF_{1,0};  \\
 \gd_f&:=& \gd_{f(1)}\t \cdots \t \gd_{f(n)}\;\;
 {\rm where}\;\; \gd_0:=F_{1,0}, \;\; \gd_1:= D_1;  \\
 \gd_f'&:=& \gd_{f(1)}'\t \cdots \t \gd_{f(n)}'\;\;
 {\rm where}\;\; \gd_0':=F_{1,0}, \;\; \gd_1':= \mD_1.   \\
\end{eqnarray*}
Note that $\gf_f\in \CJ (\mF_{n,0})_{\s , \tau}$, $\gd_f\in \CJ
(D_n)_{\s , \tau}$, and $\gd_f'\in \CJ (\mD_n)_{\s , \tau}$.

\begin{lemma}\label{a11Oct9}
\begin{enumerate}
\item The map $\CC_n\ra \CJ(\mF_{n,0})_{\s, \tau}$, $C\mapsto
\gf_C:= \sum_{f\in C}\gf_f$, is a bijection where $\gf_\emptyset
:= 0$. \item The map $\CC_n\ra \CJ(D_n)_{\s, \tau}$, $C\mapsto
\gd_C:= \sum_{f\in C}\gd_f$, is a bijection where $\gd_\emptyset
:= 0$.
 \item The map $\CC_n\ra \CJ(\mD_n)_{\s, \tau}$, $C\mapsto \gd_C':=
\sum_{f\in C}\gd_f'$, is a bijection where $\gd_\emptyset' := 0$.
\end{enumerate}
\end{lemma}

{\it Proof}. 1. It follows from $\mF_{n,0} =\bigotimes_{i=1}^n
(K+\sum_{j\in \N}Ke_{jj}(i))$, $\s_i (e_{jj}(i))= e_{j+1,
j+1}(i)$, $\tau_i (e_{jj}(i))= e_{j-1, j-1}(i)$,
$e_{kk}(i)e_{jj}(i)=\d_{jk}e_{jj}(i)$ that any ideal $\gb \in \CJ
(\mF_{n,0})_{\s , \tau}$ is a sum $\sum_{f\in C'}\gf_f$. Then $\gb
= \gf_C$ for a unique element $C\in \CC_n$ ($C$ is the set of all
the maximal elements of $C'$, it does not depend on $C'$), and so
the map $C\mapsto \gf_C$ is a bijection.

2. Similarly, it follows from $D_n=\bigotimes_{i=1}^n
(K[H_i]+\sum_{j\in \N}Ke_{jj}(i))$, $\tau _i(H_i) = H_i+1$ (hence
$K[H_i]=\bigcup_{s\geq 1}\ker ( \tau_i-1)^s$)  and the actions of
the endomorphisms $\s_i$, $\tau_i$ on the matrix units $e_{jj}(i)$
that any ideal $\gb\in \CJ (D_n)_{\s , \tau}$ is a sum $\sum_{f\in
C'}\gd_f$. Then   $\gb = \gd_C$ for a unique element $C\in \CC_n$
($C$ is the set of all the maximal elements of $C'$, it does not
depend on $C'$), and so the map $C\mapsto \gd_C$ is a bijection.

3. Statement 3 follows from statement 2 since the commutative
algebra $\mD_n$ is a localization of the commutative algebra $D_n$
at the monoid generated by the set $\{ (H_i-j)_1\, | \, i=1,
\ldots , n; 0\neq j\in \N\}$ of nonzero divisors and $\gd_C'=
\mD_n \gd_C$ for all $C\in \CC_n$. $\Box $

\begin{corollary}\label{b11Oct9}
\begin{enumerate}
\item The restriction map $\CJ (D_n)_{\s , \tau }\ra \CJ
(\mF_{n,0})_{\s , \tau}$, $\gb \mapsto \gb^r:= \gb \cap
\mF_{n,0}$, is an isomorphism (i.e. $(\gb_1*\gb_2)^r= \gb_1^r
*\gb_2^r$ for $*\in \{ +, \cdot , \cap\}$) and its inverse is the
extension map $\gc\mapsto \gc^e:= D_n \gc$.
 \item The restriction map $\CJ (\mD_n)_{\s , \tau }\ra \CJ
(\mF_{n,0})_{\s , \tau}$, $\gb \mapsto \gb^r:= \gb \cap
\mF_{n,0}$, is an isomorphism (i.e. $(\gb_1*\gb_2)^r= \gb_1^r
*\gb_2^r$ for $*\in \{ +, \cdot , \cap\}$) and its inverse is the
extension map $\gc\mapsto \gc^e:= \mD_n \gc$.
\end{enumerate}
\end{corollary}

{\it Proof}. 1. Statement 1 follows from Corollary
\ref{a11Oct9}.(1,2) and the fact that $\gd^r_C= \gf_C$ for all
$C\in \CC_n$.

2. Statement 2 follows from Corollary \ref{a11Oct9}.(1,3) and the
fact that $(\gd_C')^r= \gf_C$ for all $C\in \CC_n$.  $\Box $

$\noindent $

{\it Proof of Theorem \ref{10Oct9}}. 1. By Theorem
\ref{A10Oct9}.(1), the restriction map $\CJ (\mI_n) \ra \CJ
(D_n)_{\s , \tau}$ is an isomorphism and its inverse map is the
extension map. By Corollary \ref{b11Oct9}.(1), the restriction map
$\CJ (D_n)_{\s , \tau}\ra \CJ (\mF_n)_{\s , \tau}$ is an
isomorphism and its inverse map is the extension map. Now,
statement 1 is obvious.

2. Similarly, by Theorem \ref{A10Oct9}.(2), the restriction map
$\CJ (\mA_n) \ra \CJ (\mD_n)_{\s , \tau}$ is an isomorphism and
its inverse map is the extension map. By Corollary
\ref{b11Oct9}.(2), the restriction map $\CJ (\mD_n)_{\s , \tau}\ra
\CJ (\mF_{n,0})_{\s , \tau}$ is an isomorphism and its inverse map
is the extension map. Now, statement 2 is obvious.  $\Box $

\begin{theorem}\label{16Nov9}
Let $\Id (\CI_n)$ be the set of all the idempotent ideals of the
algebra $\CI_n$. Then
\begin{enumerate}
\item the restriction map $\CI (\mI_n) \ra \Id (\CI_n)$, $\ga
\mapsto \ga^e:= \ga \cap \CI_n$ is a bijection such that
$(\ga_1*\ga_2)^r= \ga_1^r *\ga_2^r$ for $*\in \{ +, \cdot ,
\cap\}$, and its inverse is the extension map $\gb \mapsto \gb^e:=
\mI_n \gb \mI_n$. \item The restriction map $ \Id (\CI_n) \ra \CJ
(\mF_{n,0})_{\s , \tau}$, $\gb \mapsto \gb^r:= \gb \cap
\mF_{n,0}$, is a bijection such that $(\gb_1*\gb_2)^r= \gb_1^r
*\gb_2^r$ for $*\in \{ +, \cdot , \cap\}$, and its inverse  is the
extension map $\gc \mapsto \gc^e:= \CI_n \gc \CI_n$.
\end{enumerate}
\end{theorem}

{\it Proof}. 1. Statement 1 follows from Theorem \ref{10Oct9}.(1)
and statement 2.

2. Statement 2 follows at once from a classification of the
idempotent ideals of the algebra $\mS_n\simeq \CI_n$ (Theorem 7.2,
\cite{shrekalg}). $\Box $


\section{The Noetherian factor algebra of the algebra $\mI_n$}\label{NFAAI}

The aim of this section is to show that the factor algebra $\mI_n/
\ga_n$ of the algebra $\mI_n$ at its maximal ideal $\ga_n=
\gp_1+\cdots +\gp_n$ is the only Noetherian factor algebra of the
algebra $\mI_n$ (Proposition \ref{A17Oct9}).

{\bf The factor algebra $\mI_n/ \ga_n$}. Recall that the Weyl
algebra $A_n$ is the  generalized Weyl algebra  $\CP_n
((\s_1,...,\s_n), (H_1, \ldots , H_n))$.  Denote by $S_n$ the
multiplicative submonoid of $\CP_n$ generated by the elements
$H_i+j$ where  $i=1, \ldots , n$ and $j\in \mathbb{Z}$. It follows
from the above presentation of the Weyl algebra $A_n$ as a GWA
that $S_n$ is an Ore set in $A_n$, and, using the
$\mathbb{Z}^n$-grading of $A_n$, that the (two-sided) localization
$\CA_n:=S_n^{-1}A_n$ of the Weyl algebra $A_n$ at $S_n$ is the
{\em skew Laurent polynomial ring} 
\begin{equation}\label{Anskewlaurent}
\CA_n=S_n^{-1}\CP_n [x_1^{\pm 1}, \ldots ,x_n^{\pm 1};
\s_1,...,\s_n]
\end{equation}
with coefficients in the algebra
$$ \CL_n:=S_n^{-1}\CP_n=K[H_1^{\pm 1}, (H_1 \pm 1)^{-1}, (H_1 \pm 2)^{-1}, \ldots ,
H_n^{\pm 1}, (H_n \pm 1)^{-1}, (H_n \pm 2)^{-1}, \ldots ],$$ which
is the localization of $\CP_n$ at $S_n$. We identify the Weyl
algebra $A_n$ with its image in the algebra $\CA_n$ via the
monomorphism,
$$A_n\ra \CA_n, \;\; x_i\mapsto x_i,\;\; \der_i\mapsto  H_ix_i^{-1}, \;\;
i=1, \ldots , n.$$ Let $k_n$ be the $n$'th {\em Weyl skew field},
that is the full ring of quotients of the $n$'th Weyl algebra
$A_n$ (it exists by  Goldie's Theorem  since $A_n$ is a Noetherian
domain). Then the algebra $\CA_n$ is a $K$-subalgebra of $k_n$
generated by the elements $x_i$, $x_i^{-1}$, $H_i$ and $H_i^{-1}$,
$i=1, \ldots , n$ since, for all  $j\in \N$, 
\begin{equation}\label{Hijxi}
 (H_i\mp j)^{-1}=x_i^{\pm j}H_i^{-1}x_i^{\mp j}, \;\; i=1, \ldots
, n.
\end{equation}
Clearly, $\CA_n \simeq \CA_1 \t \cdots \t \CA_1$ ($n$ times).

Recall that the algebra $\mI_n$ is a subalgebra of $\mA_n$ and the
extension $\ga_n^e$ of the maximal ideal $\ga_n$ of the algebra
$\mI_n$ is the maximal ideal of the algebra $\mA_n$. By (22) of
\cite{Bav-Jacalg}, there is the algebra isomorphism (where $\oa :=
a+\ga_n^e$):
$$\mA_n/ \ga_n^e\ra \CA_n, \;\; \overline{x_i}\mapsto x_i, \;\;
\overline{\der}_i\mapsto H_ix_i^{-1}, \;\; \overline{H_i^{\pm
1}}\mapsto H_i^{\pm 1}, \;\; i=1, \ldots, n. $$ Since
$\ga_n^{er}=\ga_n$ (Theorem \ref{7Oct9}), the algebra $B_n:=
\mI_n/ \ga_n$ is a subalgebra of the algebra $\mA_n/ \ga_n^e$, and
so there is the algebra monomorphism (where  $\oa := a+\ga_n^e$):
$$ B_n\ra \CA_n, \;\; \overline{x_i}\mapsto x_i, \;\;
\overline{\der}_i\mapsto H_ix_i^{-1}, \;\;
\overline{\int_i}\mapsto x_iH_i^{- 1}, \;\; \overline{H_i}\mapsto
H_i, \;\; i=1, \ldots, n. $$ It follows that there is the algebra
isomorphism: $$B_n\ra \bigotimes_{i=1}^n K[H_i][\der_i,
\der_i^{-1}; \tau_i]=\CP_n [ \der_1^{\pm 1}, \ldots, \der_n^{\pm
1}; \tau_1, \ldots , \tau_n],$$ the RHS is the skew Laurent
polynomial algebra with coefficients in the polynomial algebra
$\CP_n= K[H_1, \ldots , H_n]$ where $\tau_i(H_j) = H_j+\d_{ij}$.
It is a standard fact that 
\begin{equation}\label{BnL1}
B_n=(A_n)_{\der_1, \ldots , \der_n}
\end{equation}
where $(A_n)_{\der_1, \ldots , \der_n}$ is the localization of the
Weyl algebra $A_n$ at the Ore subset of $A_n$ which is the
submonoid of $A_n$  generated by the elements $\der_1, \ldots ,
\der_n$. Note that $(A_n)_{\der_1, \ldots , \der_n}\simeq
(A_n)_{x_1, \ldots , x_n}$. It is well-known that the algebra
$B_n$ is a simple, Noetherian, finitely generated algebra of
Gelfand-Kirillov dimension $2n$ and $\lgldim (B_n) = \rgldim (B_n)
= n$.

\begin{proposition}\label{A17Oct9}
Let $\ga$ be an ideal of the algebra $\mI_n$ such that $\ga \neq
\mI_n$. The following statements are equivalent.
\begin{enumerate}
\item  The factor algebra $\mI_n/\ga$ is a left Noetherian
algebra. \item The factor algebra $\mI_n/\ga$ is a right
Noetherian algebra. \item The factor algebra $\mI_n/\ga$ is a
Noetherian algebra. \item $\ga = \ga_n$.
\end{enumerate}
\end{proposition}

{\it Proof}. Note that the algebra $B_n=\mI_n/\ga_n$ is a
Noetherian algebra as a two-sided localization of the Noetherian
algebra $A_n$. Suppose that $\ga \neq \ga_n$. Fix $\gp \in \Min
(\ga )$. Then $\gp = \gp_I:=\sum_{i\in I} \gp_i$ for a non-empty
subset $I$ of the set $\{ 1, \ldots , n\}$ with $m:= |I|<n$
(Corollary \ref{b10Oct9}.(10) and Corollary \ref{c10Oct9}). The
factor algebra $\mI_n / \gp\simeq B_m\t \mI_{n-m}$ is neither left
nor right Noetherian since the algebra $\mI_{n-m}$  is so. The
algebra $\mI_n/\gp$ is a factor algebra of the algebra $\mI_n /
\ga$. Then the algebra $\mI_n / \ga$ is neither left nor right
Noetherian. Now, the proposition is obvious.  $\Box $

\begin{lemma}\label{a18Oct9}
Let $\ga$ be an ideal of the algebra $\mI_n$ distinct from
$\mI_n$. Then $\GK (\mI_n / \ga ) = 2n$.
\end{lemma}

{\it Proof}. It is well-known that $\GK (B_n) = 2n$. Now, $2n =
\GK (\mI_n) \geq \GK ( \mI_n / \ga ) \geq \GK (\mI_n / \ga_n ) =
\GK (B_n) = 2n$. Therefore, $\GK ( \mI_n / \ga ) = 2n$. $\Box $

$\noindent $


\section{The group of units of the algebra $\mI_n$ and its centre}\label{GUAIC}

In this section, the group $\mI_n^*$ of units of the algebra
$\mI_n$ is described (Theorem \ref{17Oct9}.(1)) and its centre is
found (Theorem \ref{17Oct9}.(2)). It is proved that the algebra
$\mI_n$ is central  (Lemma \ref{c12Oct9}.(2))  and self-dual.

{\bf The involution $*$ on the algebra $\mI_n$}. Using the
defining relations in Proposition \ref{a5Oct9}.(1), we see that
the algebra $\mI_n$ admits the involution: 
\begin{equation}\label{*invIn}
*: \mI_n\ra \mI_n, \;\; \der_i \mapsto  \int_i, \;\; \int_i\mapsto
\der_i, \;\; H_i\mapsto H_i, \;\; i=1, \ldots , n,
\end{equation}
i.e. it is a $K$-algebra {\em anti-isomorphism} $((ab)^*= b^*a^*)$
such that $*\circ *= \id_{\mI_n}$. Therefore, the algebra $\mI_n$
is {\em self-dual}, i.e. is isomorphic to its {\em opposite}
algebra $\mI_n^{op}$. As a result, the left and the right
properties of the algebra $\mI_n$ are the same.  For all elements
$\alpha , \beta \in \N^n$, 
\begin{equation}\label{eab*}
e_{\alpha\beta}^*= e_{\beta \alpha}.
\end{equation}
An element $a\in \mI_n$ is called {\em hermitian} if $a^*=a$.

\begin{lemma}\label{b12Oct9}
\begin{enumerate}
\item $\ga^* = \ga$ for all ideals $\ga$ of the algebra
$\mI_n$.\item $(\mI_{n,\alpha})^*= \mI_{n, -\alpha}$ for all
$\alpha \in \Z^n$. \item The set ${\rm Fix}_{\mI_n}(*)=\{ a\in
\mI_n\, | \, a^*=a\}$ of all the hermitian elements of the algebra
$\mI_n$ is the commutative subalgebra $D_n$ of the algebra
$\mI_n$.
\end{enumerate}
\end{lemma}

{\it Proof}. 1. By (\ref{eab*}), $\gp_i^*= \gp_i$ for all $i=1,
\ldots , n$ (see Corollary \ref{b10Oct9}.(2)). By Corollary
\ref{b10Oct9}.(4,9), $\ga^* = \ga$.

2. Note that $D_n^* = D_n$ and $v_{\alpha}^* = v_{-\alpha}$. By
Proposition \ref{a5Oct9}.(2), $(\mI_{n,\alpha})^* = (v_\alpha
D_n)^* = D_nv_{-\alpha}= \mI_{n, -\alpha}$.

3. By statement 2, ${\rm Fix}_{\mI_n}(*)\subseteq D_n$. The
opposite inclusion is obvious. Therefore, ${\rm
Fix}_{\mI_n}(*)=D_n$. $\Box $

$\noindent $

The involution $*$ of the algebra $\mI_n$ respects the maximal
ideal $\ga_n$ ($\ga_n^* = \ga_n$). Therefore, the factor algebra
$B_n= \mI_n / \ga_n$ inherits the involution $*$: $\der_i^* =
\der_i^{-1}$, $x_i^* = x_i +\der_i^{-1}$, $H_i^* = H_i$ for $i=1,
\ldots , n$ (since $\der_i^* = \overline{\int_i}=\der_i^{-1}$ and
$x_i^* =(\der_i H_i)^*= H_i\der_i^{-1} = \der_i x_i \der_i^{-1} =
x_i +\der_i^{-1}$).

The involution $*$ of the algebra $\mI_n$  can be extended to an
involution of the algebra $\mA_n$ by setting
$$ x_i^*=H_i\der_i, \;\; \der_i^* = \int_i, \;\; (H_i^{\pm 1})^*=
H_i^{\pm 1}, \;\; i=1, \ldots , n.$$ This can be checked using the
defining relations coming from the presentation of the algebra
$\mA_n$ as a GWA. Note that $y_i^* = (H_i^{-1}\der_i)^* = \int_i
H_i^{-1} = x_iH_i^{-2}$, $A_n^* \not\subseteq A_n$, $\mS_n^*
\not\subseteq \mS_n$, but $\CI_n^* = \CI_n$ where $\CI_n$ is the
algebra of integro-differential operators with constant
coefficients.

For a subset $S$ of a ring $R$, the sets $\lann_R(S):= \{ r\in R\,
| \, rS=0 \}$ and $\rann_R(S):= \{ r\in R\, | \, Sr=0\}$ are
called the {\em left} and the {\em right annihilators} of the set
$S$ in $R$.
 Using the fact that the algebra $\mI_n$ is a GWA and its
 $\Z^n$-grading, we see that
\begin{equation}\label{laI}
\lann_{\mI_n}(\int_i) = \bigoplus_{k\in \N}Ke_{k0}(i)\bigotimes
\bigotimes_{i\neq j}\mI_1(j), \;\; \rann_{\mI_n}(\int_i) =0.
\end{equation}
\begin{equation}\label{lad}
\rann_{\mI_n}(\der_i) = \bigoplus_{k\in \N}Ke_{0k}(i)\bigotimes
\bigotimes_{i\neq j}\mI_1(j), \;\; \lann_{\mI_n}(\der_i) =0.
\end{equation}
Recall that a submodule of a module that intersects non-trivially
each nonzero submodule of the module is called an {\em essential}
submodule.

\begin{lemma}\label{a17Oct9}
\begin{enumerate}
\item For all nonzero ideals $\ga$ of the algebra $\mI_n$,
$\lann_{\mI_n}(\ga ) = \rann_{\mI_n}(\ga) =0$. \item Each nonzero
ideal of the algebra $\mI_n$ is an essential left and right
submodule of $\mI_n$.
\end{enumerate}
\end{lemma}

{\it Proof}.  The algebra $\mI_n$ is self-dual, so it suffices to
prove only, say, the left versions of the statements.

1. Suppose that $\gb := \lann_{\mI_n}(\ga )\neq 0$, we seek a
contradiction. By Corollary \ref{b10Oct9}.(8), the nonzero ideals
$\ga$ and $\gb$ contain the ideal $F_n$. Then $0=\gb \ga \supseteq
F_n^2= F_n\neq 0$, a contradiction. Therefore, $\gb =0$.

2. Let $I$ be a nonzero left ideal of the algebra $\mI_n$. By
statement 1, $0\neq F_n I \subseteq F_n \cap I$. Therefore, $F_n$
is an essential left submodule of the algebra $\mI_n$. Then so are
all the nonzero ideals of the algebra $\mI_n$ since $F_n$ is the
least nonzero ideal of the algebra $\mI_n$.  $\Box $

\begin{corollary}\label{a24Oct9}
Let $A$ be a $K$-algebra. Then the algebra $\mI_n\t A$ is a prime
algebra iff the algebra $A$ is so.
\end{corollary}

{\it Proof}. It is obvious that if the algebra $A$ is not prime
($\ga \gb=0$ for some nonzero ideals $\ga $ and $\gb$  of $A$)
then the algebra $\mI_n\t A$ is neither (since $\mI_n \t \ga \cdot
\mI_n\t \gb =0$).

It suffices to show that if the algebra $A$ is prime then so is
the algebra $\mI_n\t A$. Let $\gc$ be a nonzero ideal of the
algebra $\mI_n\t A$. Then $F_n\gc\neq 0$, by Lemma
\ref{a17Oct9}.(1). Note that $F_n\gc \subseteq \gc$. Let
$u=E_{\alpha \beta} \t a+\cdots + E_{\s \rho}\t a'$ be a nonzero
element of $F_n\gc$ where $E_{\alpha\beta}, \ldots , E_{\s\rho}$
are distinct matrix units;  $a, \ldots , a'\in A$, and $a\neq 0$.
Then $0\neq E_{\alpha \beta} \t a= E_{\alpha\alpha}
uE_{\beta\beta} \in \ga$, and so $F_n\t AaA\subseteq \gc$.
 Let $\gd$ be a nonzero ideal of the algebra $\mI_n\t A$. Then
 $F_n\t AbA\subseteq \gd$ for some nonzero element $b\in
 A$. Then
$$\gc\gd \supseteq F_n\t AaA\, \cdot \, F_n\t AbA= F_n\t (
AaA\cdot AbA)\neq 0$$ since $F_n^2= F_n$ and $AaA\cdot AbA\neq 0$
($A$ is a prime algebra). Therefore, $\mI_n\t A$ is a prime
algebra. $\Box $

{\bf The centre of the algebra $\mI_n$}.  For an algebra $A$ and
its subset $S$, the subalgebra of $A$, $\Cen_A(S):=\{ a\in A\, |
\, as=sa$ for all $s\in S\}$, is called the {\em centralizer} of
$S$ in $A$. The next lemma shows that the algebra $\mI_n$ is a
{\em central} algebra, i.e. its centre $Z(\mI_n)$ is $K$.

\begin{lemma}\label{c12Oct9}
\begin{enumerate}
\item $\Cen_{\mI_n}(\mF_{n,0})= \Cen_{\mI_n}(D_n) = D_n$. \item
 The centre of the algebra $\mI_n$ is $K$.
 \item $\Cen_{\mI_n}(\CI_n) = K$.
\end{enumerate}
\end{lemma}

{\it Proof}. 1. Since $\mF_{n,0}\subset D_n$ and $D_n$ is a
commutative algebra, we have the inclusions $D_n \subseteq
\Cen_{\mI_n}(D_n) \subseteq \Cen_{\mI_n}(\mF_{n,0})$. It remains
to show that the inclusion $C:= \Cen_{\mI_n}(\mF_{n,0})\subseteq
D_n$ holds. Recall that the algebra $\mI_n =\bigoplus_{\alpha \in
\Z^n} \mI_{n, \alpha}$ is a $\Z^n$-graded algebra with
$\mF_{n,0}\subset D_n =\mI_{n,0}$. Therefore, $C$ is a homogeneous
subalgebra of $\mI_n$, i.e. $C=\bigoplus_{\alpha \in \Z^n}
C_{\alpha}$ where $C_{\alpha } := C\cap \mI_{n,\alpha}$. We have
to show that $C_\alpha =0$ for all $\alpha \neq 0$. Let $c\in
C_\alpha$ for some $\alpha\neq 0$. Then $c=v_{\alpha,
+}dv_{\alpha, -}$ for some element $d\in D_n$ (the elements
$v_{\alpha, +}$ and $v_{\alpha, -}$ are defined in Theorem
\ref{A10Oct9}.(1)). For all elements  $E_{\beta \beta}\in
\mF_{n,0}$ where $\beta \in \N^n$,
\begin{eqnarray*}
 cE_{\beta\beta}&=&v_{\alpha, +}d\tau_{\alpha , -}
 (E_{\beta\beta})v_{\alpha, -}=
 v_{\alpha, +}dE_{\beta-\alpha_-, \beta-\alpha_-}v_{\alpha, -},  \\
 E_{\beta\beta}c&=&v_{\alpha, +}\tau_{\alpha , +}
 (E_{\beta\beta})dv_{\alpha, -}=
 v_{\alpha, +}E_{\beta-\alpha_+, \beta-\alpha_+}dv_{\alpha, -},  \\
\end{eqnarray*}
where $\tau_{\alpha, -}:=\prod_{\alpha_i<0}\tau_i$, $\tau_{\alpha,
+}:=\prod_{\alpha_i>0}\tau_i$,
$\alpha_-:=-\sum_{\alpha_i<0}\alpha_ie_i$ and
$\alpha_+:=\sum_{\alpha_i>0} \alpha_ie_i$  ($E_{st}=0$ if either
$s\not\in\N^n$ or $t\not\in \N^n$). Since $cE_{\beta\beta} =
E_{\beta \beta}c$ and the map $a\mapsto v_{\alpha, +}av_{\alpha,
-}$ is injective (its left inverse is the map $a\mapsto u_{\alpha,
-}au_{\alpha , +}$, see (\ref{uua})), we have the equality
$E_{\beta - \alpha_+, \beta -\alpha_+}d=E_{\beta - \alpha-, \beta
-\alpha_-}d$ for each $\beta \in \N^n$. Since $\bigoplus_{\g \in
\N^n} KE_{\g\g}$ is the direct sum of ideals of the algebra $D_n$,
it follows that $E_{\g\g}d=0$  for all elements $\g \in \N^n$.
Then it is not difficult to show that $d=0$ (using the fact that
each polynomial of $K[H_1, \ldots , H_n]$ is uniquely determined
by its values on the set $\N^n$).

 2. By statement 1, the centre $Z$ of the algebra $\mI_n$ is a
 subalgebra of $D_n$. Let $d\in Z$. For all elements $i=1, \ldots
 , n$, $0=dx_i-x_id= x_i(\tau_i(d) - d)$. Since $\mI_n\subseteq
 \mA_n$, we see that $0=y_ix_i (\tau_i(d) - d)= \tau_i(d) -d$, and
 so $d\in \bigcap_{i=1}^n \ker_{D_n}(\tau_i-1)=K$. Therefore, $Z=K$.

3. By (\ref{eijEij}), $\mF_{n,0}\subseteq \CI_n$. This implies
that $C:=\Cen_{\mI_n}(\CI_n)\subseteq
\Cen_{\mI_n}(\mF_{n,0})=D_n$, by statement 1. Let $d\in C$. Then
$$ 0=\der_i\cdot 0=\der_i(d\int_i-\int_id) =
\der_i\int_i(\tau_i(d)-d)=\tau_i(d)-d\;\; {\rm for\; all}\;\; i=1,
\ldots , n,$$ where $\tau_i(a) = \der_ia\int_i$. Hence $d\in
\bigcap_{i=1}^n \ker_{D_n}(\tau_i-1)=K$, and so $C=K$.   $\Box $

\begin{lemma}\label{a22Nov9}
Let $C=P_n$, $K[\der_1, \ldots, \der_n]$, $K[\int_1, \ldots ,
\int_n]$ or $D_n$. Then $\Cen_{\mI_n}(C) = C$ and $C$ is a maximal
 commutative subalgebra of the algebra $\mI_n$.
\end{lemma}

{\it Proof}. The first statement, $\Cen_{\mI_n}(C) = C$, follows
from the fact that the algebra $\mI_n$ is $\Z^n$-graded and the
canonical generators of the algebra $C$ are homogeneous elements
of the algebra $\mI_n$ (we leave this as an exercise for the
reader). Then $C$ is a maximal
 commutative subalgebra of the algebra $\mI_n$ since $\Cen_{\mI_n}(C) = C$
 and $C$ is a commutative algebra.
$\Box $


{\bf The group $\mI_n^*$ of units of the algebra $\mI_n$ and its
centre}.  The group $\mA_1^*$ of units of the algebra $\mA_1$
contains the following infinite discrete subgroup Theorem 4.2,
\cite{Bav-Jacalg}:
\begin{equation}\label{defH1}
\CH := \{ \prod_{i\geq 0} (H+i)^{n_i}\cdot \prod_{i\geq
1}(H-i)^{n_{-i}}_1\, | \, (n_i)_{i\in \Z}\in \Z^{(\Z )}\}\simeq
\Z^{(\Z )}.
\end{equation}
 For each
tensor multiple $\mA_1(i)$ of the algebra
 $\mA_n=\bigotimes_{i=1}^n\mA_1(i)$, let $\CH_1 (i)$ be the corresponding
group $\CH$. Their (direct) product 
\begin{equation}\label{defH2}
 \CH_n:= \CH_1 (1) \cdots \CH_1
(n)= \prod_{i=1}^n\CH_1 (i)
\end{equation}
 is a (discrete) subgroup of the group
$\mA_n^*$  of units of the algebra $\mA_n$, and $\CH_n\simeq
\CH^n\simeq (\Z^n)^{(\Z )}$. Note that $\mA_n^* = K^* \times
(\CH_n\ltimes (1+\ga_n^e)^*)$ and $Z(\mA_n^*) = K^*$ (Theorem 4.4,
\cite{Bav-Jacalg}). A similar result holds for the group $\mI_n^*$
of the algebra $\mI_n$ (Theorem \ref{17Oct9}). Since $\ga_n$ is an
ideal of the algebra $\mI_n$, the intersection $(1+\ga_n)^*:=
\mI_n^*\cap (1+\ga_n)$ is a subgroup of the group $\mI_n^*$ of
units of the algebra $\mI_n$.

\begin{theorem}\label{17Oct9}
\begin{enumerate}
\item Let $\mF_n := \bigoplus_{i=1}^n (K+F(i))$. Then $$\mI_n^*=
K^* \times (1+\ga_n )^*\;\; {\rm  and }\;\; \mI_n^* \supseteq
(1+\mF_n\cap \ga_n )^* \simeq \underbrace{\GL_\infty
(K)\ltimes\cdots \ltimes \GL_\infty (K)}_{2^n-1 \;\; {\rm
times}}.$$ \item The centre of the group $\mI_n^*$ is $K^*$.
\end{enumerate}
\end{theorem}

{\it Proof}. 1. The commutative diagram of  algebra homomorphisms
$$\xymatrix{\mI_n\ar[d]\ar[r] & \mA_n\ar[d]\\
B_n\ar[r] &\CA_n}$$ yields the commutative diagram of group
homomorphisms $$\xymatrix{\mI_n^*\ar[d]\ar[r] & \mA_n^*\ar[d]\\
B_n^*\ar[r] &\CA_n^*.}$$ Since $B_n^* =\bigcup_{\alpha \in \Z^n}
K^* \der^\alpha$ and $\mA_n^* = K^* \times (\CH_n\ltimes
(1+\ga_n^e)^*)$, we see that
$$
 K^*\times (1+\ga_n )^*\subseteq \mI_n^* \subseteq  \mI_n\cap \mA_n^*
 = K^* \times (\mI_n\cap (1+\ga_n^e)^*)
 = K^* \times (1+\ga_n^{er})^*= K^* \times (1+\ga_n )^*$$
 since $\ga_n^{er} = \ga_n$ (Theorem \ref{7Oct9}). Therefore,
 $\mI_n^* = K^*\times (1+\ga_n)^*$.

 Since $\mF_n\subset \mI_n\subset \mA_n$, it is obvious that
 $$\mI_n^* \supseteq
(1+\mF_n\cap \ga_n )^*  = (1+\mF_n\cap \ga_n^e)^*\simeq
\underbrace{\GL_\infty (K)\ltimes\cdots \ltimes \GL_\infty
(K)}_{2^n-1 \;\; {\rm times}}.$$ The isomorphism is established in
Corollary 7.3, \cite{jacaut}.

 2. Let $S$ be the set of elements
of the type $1+\prod_{i\in I}e_{s_i s_i}(i)$ where $\emptyset \neq
I\subseteq \{ 1, \ldots , n\}$. Then $S\subseteq \mI_n^*$ and
$\Cen_{\mI_n} (S) = \Cen_{\mI_n}(\mF_{n,0}) = D_n$, by Lemma
\ref{c12Oct9}.(1). Therefore, $\Cen_{\mI_n^*} (S)= \Cen_{\mI_n}
(S)\cap \mI_n^* = D_n \cap \mI_n^* = D_n^* = \mF_{n,0}^*$. We see
that  $\Cen_{\mF_{n, 0}} (S)=K$. Therefore, the centre of the
group $\mI_n^*$ is $K^*$.  $\Box $

{\bf The group of units $(1+F)^*$ and $\mI_1^*$}. Recall that the
 algebra (without 1) $F=\bigoplus_{i,j\in \N} Ke_{ij}$ is the union
$M_\infty (K) := \bigcup_{d\geq 1}M_d(K)= \varinjlim M_d(K)$ of
the matrix algebras $M_d(K):= \bigoplus_{1\leq i,j\leq
d-1}Ke_{ij}$, i.e. $F= M_\infty (K)$.  For each $d\geq 1$,
consider the (usual) determinant $\det_d=\det : 1+M_d(K)\ra K$,
$u\mapsto \det (u)$. These determinants determine the (global)
{\em determinant}, 
\begin{equation}\label{gldet}
\det : 1+M_\infty (K)= 1+F\ra K, \;\; u\mapsto \det (u),
\end{equation}
where $\det (u)$ is the common value of all the determinants
$\det_d(u)$, $d\gg 1$. The (global) determinant has usual
properties of the determinant. In particular, for all $u,v\in
1+M_\infty (K)$, $\det (uv) = \det (u) \cdot \det (v)$. It follows
from this equality and the  Cramer's formula for the inverse of a
matrix that the group $\GL_\infty (K):= (1+M_\infty (K))^*$ of
units of the monoid $1+M_\infty (K)$ is equal to 
\begin{equation}\label{GLiK}
\GL_\infty (K) = \{ u\in 1+M_\infty (K) \, | \, \det (u) \neq 0\}.
\end{equation}
Therefore, 
\begin{equation}\label{1GLiK}
(1+F)^* = \{ u\in 1+F \, | \, \det (u) \neq 0\}=\GL_\infty (K).
\end{equation}

\begin{corollary}\label{b17Oct9}
$\mI_1^* = K^* \times (1+F)^* = K^* \times \GL_\infty (K)$, i.e.
$\mI_1^* = \{ \l (1+f)\, | \, \det (1+f) \neq 0, \l \in K^*, f\in
F \}$. The elements $\l \in K^*$, $1+\mu e_{ij}$ where $\mu \in K$
and $i\neq j$, and $1+\g e_{00}$ where $\g \in K\backslash \{
-1\}$ are generators for the group $\mI_n^*$.
\end{corollary}


\section{The weak and the global dimensions of the algebra
$\mI_n$}\label{WGDAI}

In this section, we prove that the weak  dimension of the algebra
$\mI_n$ and of all its prime factor algebras is  $n$ (Theorem
\ref{18Oct9}). An analogue of  Hilbert's Syzygy Theorem  is
established for the algebra $\mI_n$ and for all its prime factor
algebras (Theorem \ref{19Oct9}).

{\bf The  weak dimension of the algebra $\mI_n$}. Let $S$ be a
non-empty multiplicatively closed subset of a ring $R$, and let
$\ass (S):= \{ r\in R\, | \, sr =0$ for some $s\in S\}$. Then a
{\em left quotient ring} of $R$ with respect to $S$ is a ring $Q$
together with a homomorphism $\v : R\ra Q$ such that

(i) for all $s\in S$, $\v (s)$ is a unit in $Q$;

(ii) for all $q\in Q$, $q= \v (s)^{-1} \v (r)$ for some $r\in R$
and $s\in S$, and

(iii) $\ker (\v ) = \ass (S)$.

If there exists a left quotient ring $Q$ of $R$ with respect to
$S$ then it is unique up to isomorphism, and it is denoted
$S^{-1}R$. It is also said that the ring $Q$ is the {\em left
localization} of the ring $R$ at $S$.

{\it Example 1}. Let $S:= S_\der:= \{ \der^i, i\geq 0\}$ and
$R=\mI_1$. Then  $\ass (S) =F$, $\mI_1/\ass (S)= B_1$ and the
conditions (i)-(iii) hold where $Q=B_1$. This means that the ring
$B_1=\mI_1/F$ is the left quotient ring of $\mI_1$ at $S$, i.e.
 $B_1\simeq S^{-1}_\der \mI_1$.

{\it Example 2}. Let $S:=S_{\der_1, \ldots , \der_n}:= \{
\der^\alpha , \alpha \in \N^n\}$ and $R=\mI_n$. Then $\ass
(S_{\der_1, \ldots , \der_n}) = \ga_n$,  $\mI_n / \ga_n = B_n$,
and 
\begin{equation}\label{Sy1yn}
S_{\der_1, \ldots , \der_n}^{-1}\mI_n\simeq B_n,
\end{equation}
i.e. $B_n$ is the left quotient ring of $\mI_n$ at $S_{\der_1,
\ldots , \der_n}$. Note that the right localization $\mI_n
S_{\der_1, \ldots , \der_n}^{-1}$ of $\mI_n$ at $S_{\der_1, \ldots
, \der_n}$ does not exist. Otherwise, we would have $S_{\der_1,
\ldots , \der_n}^{-1}\mI_n\simeq \mI_n S_{\der_1, \ldots ,
\der_n}^{-1}$ but all the elements $\der^\alpha$ are left regular,
and we would have a monomorphism $\mI_n\ra S_{\der_1, \ldots ,
\der_n}^{-1}\mI_n\simeq B_n$, which would be impossible since the
elements $\der_i$ of  the algebra $\mI_n$ are not regular. By
applying the involution $*$ to  (\ref{Sy1yn}), we see that
\begin{equation}\label{1Sy1yn}
\mI_nS^{-1}_{\int_1, \ldots , \int_n}\simeq B_n,
\end{equation}
i.e. the algebra $B_n$ is the right localization of $\mI_n$ at the
multiplicatively closed set $S_{\int_1, \ldots , \int_n}:= \{
\int^\alpha \, | \, \alpha \in \N^n\}$.

$\noindent $

Given  a ring $R$ and modules ${}_RM$ and $N_R$, we denote by $\pd
({}_RM)$ and $\pd (N_R)$ their projective dimensions. Let us
recall a result which will be used repeatedly in  proofs later.

It is obvious that $P_n\simeq A_n/\sum_{i=1}^n A_n \der_i$. A
similar result is true for the $\mI_n$-module $P_n$ (Proposition
\ref{C17Oct9}.(2)). Note that $\pd_{A_n}(P_n)=n$ but
$\pd_{\mI_n}(P_n)=0$ (Proposition \ref{C17Oct9}.(3)).

\begin{proposition}\label{C17Oct9}
\begin{enumerate}
\item $\mI_1=\mI_1\der \bigoplus \mI_1e_{00}$ and $\mI_1= \int
\mI_1\bigoplus e_{00}\mI_1$. \item ${}_{\mI_n}P_n\simeq \mI_n /
\sum_{i=1}^n \mI_n \der_i$. \item The $\mI_n$-module $P_n$ is
projective. \item $F_n =F^{\t n}$ is a left and right projective
$\mI_n$-module. \item The projective dimension of the left and
right $\mI_n$-module $\mI_n / F_n$ is $1$. \item For each element
$\alpha \in \N^n$, the $\mI_n$-module $\mI_n / \mI_n \der^\alpha$
is projective, moreover, $\mI_n / \mI_n\der^\alpha\simeq
\bigoplus_{i=1}^n(K[x_i]\bigotimes \bigotimes_{j\neq
i}\mI_1(i))^{\alpha_j}$.
\end{enumerate}
\end{proposition}

{\it Proof}. 1. Using  the equality $\int\der = 1-e_{00}$, we see
that $\mI_1= \mI_1\der +\mI_1 e_{00}$. Since $\der e_{00} =0$ and
$e_{00}^2= e_{00}$, we have $\mI_1\der \cap \mI_1e_{00} =
(\mI_1\der \cap \mI_1e_{00})e_{00} \subseteq \mI_1\der e_{00} =0$.
Therefore, $\mI_1 = \mI_1\der \bigoplus \mI_1e_{00}$. Then
applying the involution $*$ to this equality we obtain the
equality $\mI_1= \int \mI_1\bigoplus e_{00}\mI_1$.

2. Since ${}_{\mI_1}P_1\simeq \mI_1e_{00} = \bigoplus_{i\in \N}
Ke_{i0}$, $1\mapsto e_{00}$, we have ${}_{\mI_1}P_1\simeq \mI_1 /
\mI_1 \der$, by statement 1. Therefore, ${}_{\mI_n}P_n \simeq
\bigotimes_{i=1}^n P_1(i)  \simeq \bigotimes_{i=1}^n \mI_1 (i) /
\mI_1(i) \der_i \simeq \mI_n / \sum_{i=1}^n \mI_n \der_i$.

3. By statement 1, $$\mI_n =\bigotimes_{i=1}^n \mI_1(i) =
\bigotimes_{i=1}^n (\mI_1(i) \der_i \bigoplus \mI_1(i)
e_{00}(i))=\mI_n \prod_{i=1}^n e_{00}(i) \bigoplus
  (\sum_{i=1}^n \mI_n \der_i) \simeq P_n
  \bigoplus(\sum_{i=1}^n \mI_n \der_i) .$$
Therefore, $P_n$ is a projective $\mI_n$-module.

4. Note that the left $\mI_1$-module $F= \bigoplus_{i\geq 0}
\mI_1E_{ii}\simeq \bigoplus_{i\geq 0} P_1$ is projective by
statement 2. Therefore, $F_n = F^{\t n}$ is a projective left
$\mI_n$-module. Since the ideal $F_n$ is stable under the
involution $*$, $F_n^* = F_n$,  the right $\mI_n$-module $F_n$ is
projective.

5. The short exact sequence of left and right $\mI_n$-modules
$0\ra F_n \ra \mI_n\ra \mI_n / F_n\ra 0$ does not split since
$F_n$ is an essential left and right submodule of $\mI_n$ (Lemma
\ref{a17Oct9}.(2)). By statement 4, the projective dimension of
the left and right $\mI_n$-module $\mI_n / F_n$ is 1.

6. Let $\Z^n =\bigoplus_{i=1}^n \Z e_i$ where $e_1, \ldots , e_n$
is the canonical free $\Z$-basis for $\Z^n$. Let $m=|\alpha |$.
Fix a chain of elements of $\Z^n$, $\beta_0=0, \beta_1, \ldots ,
\beta_m=\alpha$ such that, for all $i$, $\beta_{i+1}=\beta_i+e_j$
for some index $j=j(i)$. Then all the factors of the chain of left
ideals
$$\mI_n\der^\alpha =\mI_n\der^{\beta_m} \subset \mI_n\der^{\beta_{m-1}}
\subset\cdots \subset \mI_n\der^{\beta_1} \subset \mI_n$$ are
projective $\mI_n$-modules since
$\mI_n\der^{\beta_i}/\mI_n\der^{\beta_{i+1}}\simeq
\mI_n/\mI_n\der_j\simeq K[x_j]\t\mI_{n-1}$ is the projective
$\mI_n$-module (statement 3). The first isomorphism is due to the
fact that the element $\der_i$ is left regular, i.e. $a\der_i =
b\der_i$ implies $a=b$ (by multiplying the equation on the right
by $\int_i$). Therefore, the $\mI_n$-module $\mI_n/
\mI_n\der^\alpha $ is projective. Moreover, $\mI_n /
\mI_n\der^\alpha\simeq \bigoplus_{i=1}^n(K[x_i]\bigotimes
\bigotimes_{j\neq i}\mI_1(i))^{\alpha_j}$. $\Box $


\begin{theorem}\label{18Oct9}
Let $\mI_{n,m}:= B_{n-m}\t \mI_m$ where $m=0, 1, \ldots , n$ and
$\mI_0 = B_0:=K$. Then $\wdim (\mI_{n,m}) =  n$ for all $m=0, 1,
\ldots , n$. In particular, $\wdim (\mI_n)  = n$.
\end{theorem}

{\it Proof}. The algebra $B_n$ is Noetherian, hence $
 n=\lgldim (S^{-1}_{\der_1, \ldots ,
\der_n}\mI_{n,m})=\wdim (B_n) \leq \wdim (\mI_{n,m})$  (Corollary
7.4.3, \cite{MR}). To finish the proof of the theorem it suffices
to show that the inequality $\wdim (\mI_{n,m}) \leq n$ holds for
all numbers $n$ and $m$. We use induction on $n$. The case $n=0$
is trivial. So, let $n\geq 1$ and we assume the inequality holds
for all $n'<n$ and all $m=0, 1, \ldots , n'$. For $n$, we use the
second induction on $m=0, 1, \ldots , n$. When $m=0$, the
inequality holds since $\mI_{n,0}=B_n$ and $\wdim (B_n) =n$.

Suppose that $m>0$ and $\wdim (\mI_{n,m'})\leq n$ for all $m'<m$.
We have to show that $\wdim (\mI_{n,m})\leq n$ or, equivalently,
$\fd_{\mI_{n,m}}(M)\leq n$ for all $\mI_{n,m}$-modules $M$ ($\fd$
 denotes  the flat dimension). Changing the order of the tensor multiples
we can write $\mI_{n,m} = \mI_1\t \mI_{n-1, m-1}$. Then $\wdim
(\mI_{n-1, m-1})\leq n-1$, by the inductive hypothesis. Recall
that $B_1=S^{-1}_\der \mI_1= \mI_1/F$ and every $\der$-torsion
$\mI_1$-module $V$ is a direct sum of several (maybe an  infinite
number of) copies of the projective simple $\mI_1$-module $K[x]$
(Proposition \ref{C17Oct9}.(6)), hence $V$ is projective, hence
$V$ is  flat. Note that $S^{-1}_\der \mI_{n,m}\simeq \mI_{n, m-1}$
and $\wdim (\mI_{n, m-1})\leq n$, by the inductive hypothesis. The
$\mI_{n,m}$-module $ {\rm tor}_\der (M) := \{ m\in M \, | \,
\der^im =0$ for some $i\}$ is the $\der$-{\em torsion} submodule
of the $\mI_{n,m}$-module $M$. There are two  short exact
sequences of $\mI_{n,m}$-modules, 
\begin{equation}\label{tordM}
0\ra {\rm tor}_\der (M) \ra M\ra \bM \ra 0,
\end{equation}
\begin{equation}\label{1tordM}
0\ra \bM \ra S_\der^{-1}M\ra M' \ra 0,
\end{equation}
where the $\mI_{n,m}$-modules ${\rm tor}_\der (M)$ and $M'$ are
$\der$-torsion, and the $\mI_n$-module $S_\der^{-1}M$ is
$\der$-torsion free. To prove that $\fd_{\mI_{n,m}}(M)\leq n$ it
suffices to show that the flat dimensions of the
$\mI_{n,m}$-modules ${\rm tor}_\der (M)$, $S_\der^{-1}M$  and $M'$
are less or equal to $n$. Indeed, then by (\ref{1tordM}),
$\fd_{\mI_{n,m}}(\bM ) \leq \max \{ \fd_{\mI_{n,m}}(S_\der^{-1}M),
\fd_{\mI_{n,m}}(M')\}\leq n$; and by (\ref{tordM}),
$\fd_{\mI_{n,m}}(M)\leq \max \{ \fd_{\mI_{n,m}}({\rm tor}_\der
(M)), \fd_{\mI_{n,m}}(\bM )\}$ $\leq n$.

The $\mI_1$-module ${\rm tor}_\der (M)$ (where $\mI_{n,m} =
\mI_1\t \mI_{n-1, m-1}$) is a direct sum of copies (may be
infinitely many) of the projective simple $\mI_1$-module $K[x]$.
Note that $\End_{\mI_1}(K[x]) \simeq \ker_{K[x]}(\der ) =K$ since
${}_{\mI_1}K[x]\simeq \mI_1/ \mI_1\der$. Using this fact and
Proposition \ref{C17Oct9}.(6), for each finitely generated
submodule $T$ of the $\mI_{n,m}$-module ${\rm tor}_\der (M)$ there
exists a  family $\{ T_i\}_{i\in I}$ of its submodules $T_i$ where
$(I, \leq )$ is a well-ordered set such that if $i,j\in I$ and
$i\leq j$ then $T_i\subseteq T_j$, $T= \bigcup_{i\in I} T_i$ and
  $T_i/\bigcup_{j<i}T_j\simeq K[x]\t \CT_i$ for some $\mI_{n-1, m-1}$-module
   $\CT_i$. Note that
   $$\fd_{\mI_{n,m}}(K[x]\t \CT_i)
   \leq  \fd_{\mI_{n-1, m-1}}(\CT_i)\leq n-1,$$
   since the $\mI_1$-module $K[x]$ is projective.
Therefore, $\fd_{\mI_{n,m}}(T)\leq n-1$. The module ${\rm
tor}_\der (M) =\bigcup_{\theta \in \Theta} T_\theta$ is the union
of its finitely generated submodules $T_\theta$, hence
$$\fd_{\mI_{n,m}}({\rm tor}_\der (M))
=\fd_{\mI_{n,m}}(\bigcup_{\theta\in \Theta}T_\theta )\leq \sup \{
\fd_{\mI_{n,m}}(T_\theta)\}_{\theta\in \Theta}= n-1.$$
 Similarly, $\fd_{\mI_{n,m}}(\bM )\leq n-1$ since the
$\mI_{n,m}$-module $\bM$ is $\der$-torsion.

It remains to show that $\fd_{\mI_{n,m}}(S_\der^{-1}(M))\leq n$.
By (\ref{Sy1yn}), the left $\mI_1$-module $B_1$ is flat, hence the
left $\mI_{n,m}$-module $B_1\t \mI_{n-1, m-1}$ is flat. Then, by
Proposition 7.2.2.(ii), \cite{MR},
$$\fd_{\mI_{n,m}}(S^{-1}_\der M) \leq \fd_{B_1\t \mI_{n-1,
m-1}}(S^{-1}_\der M) +\fd_{\mI_{n,m}}(B_1\t \mI_{n-1, m-1})\leq
\wdim (\mI_{n, m-1})\leq n.$$ The proof of the theorem is
complete. $\Box $


\begin{corollary}\label{b29Oct9}
Let $M$ be a $\der$-torsion $\mI_{n,m}$-module, i.e. $S^{-1}_\der
M=0$, where $S^{-1}_\der : \,I_{n,m}=\mI_1\t \mI_{n-1, m-1}\ra
B_1\t \mI_{n-1, m-1}=\mI_{n,m-1}$ is the localization map and
$n,m\geq 1$. Then there exists a family $\{ T_i\}_{i\in I}$ of
$\mI_{n,m}$-submodules of $M$ such that $M=\bigcup_{i\in I}T_i$,
$(I, \leq )$ is a well-ordered set such that if $i,j\in I$ and
$i\leq j$ then $T_i\subseteq T_j$,  and
  $T_i/\bigcup_{j<i}T_j\simeq K[x]\t \CT_i$ for some $\mI_{n-1, m-1}$-module
   $\CT_i$.
\end{corollary}

{\it Proof}. The $\mI_1$-module $M$ is a direct sum of (may be
infinitely many) copies  of the projective simple $\mI_1$-module
$K[x]$. Note that $\End_{\mI_1}(K[x]) \simeq \ker_{K[x]}(\der )
=K$ since ${}_{\mI_1}K[x]\simeq \mI_1/ \mI_1\der$. Using this
fact, for the $\mI_{n,m}$-module $M$ there exists a family $\{
T_i\}_{i\in I}$ of its submodules $T_i$ where $(I, \leq )$ is a
well-ordered set such that if $i,j\in I$ and $i\leq j$ then
$T_i\subseteq T_j$, $M = \bigcup_{i\in I} T_i$ and
  $T_i/\bigcup_{j<i}T_j\simeq K[x]\t \CT_i$ for some $\mI_{n-1, m-1}$-module
   $\CT_i$. $\Box $

\begin{corollary}\label{a19Oct9}
Let $A$ be a prime factor algebra of the algebra $\mI_n$. Then
 $\wdim (A) =n$.
\end{corollary}

{\it Proof}. By Corollary \ref{b10Oct9}.(10), the algebra $A$ is
isomorphic to the algebra $\mI_{n,m}$ for some $m$. Now, the
corollary follows from Theorem \ref{18Oct9}. $\Box $

The next theorem is an analogue of  Hilbert's Syzygy Theorem for
the algebra $\mI_n$ and its prime factor algebras. The {\em flat
dimension} of an $A$-module $M$ is denoted by $\fd_A(M)$.

\begin{theorem}\label{19Oct9}
Let $K$ be an algebraically closed uncountable field of
characteristic zero. Let $A$ be a prime factor algebra of $\mI_n$
(for example, $A= \mI_n$) and $B$ be a  Noetherian finitely
generated algebra over $K$. Then
 $\wdim (A\t B) = \wdim (A) +\wdim (B)  = n+\wdim (B)$.
\end{theorem}

{\it Proof}. Recall that $A\simeq \mI_{n,m}$ for some $m\in \{0,
1, \ldots , n\}$ and $ \wdim (\mI_{n,m}) = n$ (Theorem
\ref{18Oct9}). Since
$$ n+\wdim (B) = \wdim (\mI_{n,m})+\wdim(B) \leq \wdim
(\mI_{n,m}\t B), $$ it suffices to show that $\wdim (\mI_{n,m}\t
B) \leq n+\wdim (B)$ for all numbers $n$ and $m$. We use induction
on $n$. The case $n=0$ is trivial since $A=K$. So, let $n\geq 1$,
and we assume that the inequality holds for all $n'<n$ and all
$m'=0, 1, \ldots , n'$. For the number $n\geq 1$, we use the
second induction on $m=0, 1, \dots , n$. The case $m=0$, i.e.
$\mI_{n,0}=B_n$, is known, Corollary 6.3, \cite{THM} (this can
also be deduced from Proposition 9.1.12, \cite{MR}; see also
\cite{glgwa}).

So, let $m>0$ and we assume that the inequality holds for all
numbers $m'<m$. Let $M$ be an $\mI_{n,m}\t B$-module. We have to
show that $\fd_{\mI_{n,m}\t B}(M)\leq n+\wdim (B)$. We can treat
$M$ as an $\mI_{n,m}$-module. Then we have the short exact
sequences (\ref{tordM}) and (\ref{1tordM})  which are, in fact,
short exact sequence of $\mI_{n,m}\t B$-modules. To prove that
$\fd_{\mI_{n,m}\t B}(M)\leq n+\wdim (B)$ it suffices to show that
the flat dimensions of the $\mI_{n,m}\t B$-modules ${\rm tor}_\der
(M)$, $S^{-1}_\der M$ and $M'$ are less or equal to $n+\wdim (B)$,
by the same reason as in the proof of Theorem \ref{18Oct9}.
 Repeating the same argument as at the end of the proof of Theorem
\ref{18Oct9}, for each finitely generated submodule $T$ of the
$\mI_{n,m}\t B$-module ${\rm tor}_\der (M)$ (where $\mI_{n,m} =
\mI_1\t \mI_{n-1, m-1}$) there exists a  family $\{ T_i\}_{i\in
I}$ of its submodules $T_i$ where $(I, \leq )$ is a well-ordered
set such that if $i,j\in I$ and $i\leq j$ then $T_i\subseteq T_j$,
$T = \bigcup_{i\in I} T_i$ and
  $T_i/\bigcup_{j<i}T_j\simeq K[x]\t \CT_i$ for some $\mI_{n-1, m-1}\t B$-module
   $\CT_i$. Note that $\mI_{n,m}\t B= \mI_1\t \mI_{n-1, m-1}\t B$
   and
   $$\fd_{\mI_{n,m}\t B}(K[x]\t \CT_i)
   \leq \fd_{\mI_{n-1, m-1}\t B}(\CT_i)
   \leq \wdim (\mI_{n-1, m-1}\t B)=
n-1+\wdim (B),$$ since the $\mI_1$-module $K[x]$ is projective.
 Therefore, $\fd_{\mI_{n,m}\t B}(T_i)\leq n-1+\wdim (B)$. The
 $\mI_{n,m}\t B$-module ${\rm tor}_\der (M) =\bigcup_{\theta \in
 \Theta}T_\theta$ is the union of its finitely generated
 submodules $T_\theta$, hence
$$\fd_{\mI_{n,m}\t
 B}({\rm tor}_\der (M))=\fd_{\mI_{n,m}\t B}(\bigcup_{\theta \in \Theta}T_\theta )
 \leq \sup \{ \fd_{\mI_{n,m}\t B}(T_\theta )\}_{\theta \in \Theta} = n-1+\wdim (B).$$
 Similarly, $\fd_{\mI_{n,m}\t B}(M' ) \leq n-1+\wdim (B)$ since
 the $\mI_1$-module $M'$ is $\der$-torsion.

 It remains to show that $\fd_{\mI_{n,m}\t B}(S^{-1}_\der M)\leq
 n+\wdim (B)$. By (\ref{1Sy1yn}), the left $\mI_1$-module $B_1$ is
 flat, hence the left $\mI_{n,m}\t B$-module $B_1\t \mI_{n-1,
 m-1}\t B$ is flat. Then, by Proposition 7.2.2.(ii), \cite{MR},
\begin{eqnarray*}
\fd_{\mI_{n,m}\t B}(S^{-1}_\der M)&\leq &\fd_{B_1\t
\mI_{n-1,m-1}\t
 B}(S^{-1}_\der M)+\fd_{\mI_{n,m}\t B}(B_1\t \mI_{n-1, m-1}\t
 B)\\
 &\leq & \wdim (B_1\t \mI_{n-1, m-1}\t B)\leq n-1+\wdim (B_1\t B)\\
 & =&
 n-1+\wdim (B_1) +\wdim (B)= n-1+1+\wdim (B) = n+\wdim (B),
\end{eqnarray*}
since the algebra $ B$ is Noetherian and finitely generated. The
proof of the theorem is complete.  $\Box $

{\bf The global dimension of the algebra $\mI_n$}. For all
Noetherian rings, $\wdim (M) = \gldim (M)$, and this is not true
for {\em non-Noetherian} ring, in general. For many  Noetherian
rings, including the Weyl algebras $A_n$ and the algebras $B_n$,
the known proofs of finding their global dimensions are, in fact,
about their weak dimensions as localizations and  faithfully flat
extensions  are used. This fact together with the fact that the
algebras $\mI_n$ are not Noetherian are main difficulties in
finding their global dimensions.

\begin{proposition}\label{Aus-sub}
{\rm \cite{Aus-NagoyaMJ-55}} Let $M$ be a module over an algebra
$A$, $I$ a non-empty well-ordered set, $\{ M_i\}_{i\in I}$ be  a
family of submodules of $M$ such that if $i,j\in I$ and $i\leq j$
then $M_i\subseteq M_j$. If $M=\bigcup_{i\in I}M_i$ and
$\pd_A(M_i/M_{<i})\leq n$ for all $i\in I$ where $M_{<i}:=
\bigcup_{j<i}M_j$ then $\pd_A(M)\leq n$.
\end{proposition}
Let $V\subseteq U\subseteq W$ be modules. Then the factor module
$U/V$ is called a {\em sub-factor} of the module $W$. Each algebra
$\mI_{n,m}$ is self-dual, so its left and right global dimensions
coincide. Their common value is denoted by $\gldim (\mI_{n,m})$.

\begin{proposition}\label{a29Nov9}
$n\leq \gldim (\mI_{n,m})\leq n+m$ for all $n\in \N$ and $m=0,1,
\ldots , n$. In particular, $n\leq \gldim (\mI_n) \leq 2n$.
\end{proposition}

{\it Proof}. $n=\gldim (B_n) \leq \gldim (\mI_{n,m})$, by
(\ref{Sy1yn}). It remains to show that $\gldim (\mI_{n,m})\leq
n+m$. We use induction on $n$. The case $n=0$ is trivial as
$\mI_{0,0}=K$. So, let $n\geq 1$ and we assume that the inequality
holds for all $n'<n$ and all $m=0, 1, \ldots , n'$. For $n$, we
use the second induction on $m=0, 1, \ldots , n$. When $m=0$, the
inequality holds since $\mI_{n,0}=B_n$ and $\gldim (B_n)=n$.

Suppose that $m>0$ and $\gldim (\mI_{n,m'})\leq n+m'$ for all
$m'<m$. We have to show that $\gldim (\mI_{n,m})\leq n+m$, or
equivalently $\pd_{\mI_{n,m}}(M)\leq n+m$ for all
$\mI_{n,m}$-modules $M$. Changing the order of the tensor
multiples we can write $\mI_{n,m} = \mI_1\t \mI_{n-1, m-1}$. For
the $\mI_{n,m}$-module $M$ we have the short exact sequence of
$\mI_{n,m}$-modules (\ref{tordM}).  By Corollary \ref{b29Oct9},
for the $\mI_{n,m}$-module ${\rm tor}_\der (M)$ there exists a
family $\{ T_i\}_{i\in I}$ of its submodules $T_i$ where $(I, \leq
)$ is a well-ordered set such that if $i,j\in I$ and $i\leq j$
then $T_i\subseteq T_j$, ${\rm tor}_\der (M) = \bigcup_{i\in I}
T_i$ and
  $T_i/\bigcup_{j<i}T_j\simeq K[x]\t \CT_i$ for some $\mI_{n-1, m-1}$-module
   $\CT_i$. Note that $$\pd_{\mI_{n,m}}(K[x]\t \CT_i)
   \leq \pd_{\mI_1}(K[x]) + \pd_{\mI_{n-1, m-1}}(\CT_i)\leq n+m-2.$$
By Proposition \ref{Aus-sub}, $\pd_{\mI_{n,m}}({\rm tor}_\der
(M))\leq n+m-2$. Note that $\pd_{\mI_{n,m}}(\mI_{n,m-1})\leq 1$
(since $0\ra F\ra \mI_1\ra B_1\ra 0$ is a projective resolution
for the $\mI_1$-module $B_1$) and, by Proposition 7.2.2.(ii),
\cite{MR},
$$\pd_{\mI_{n,m}}(\bM )\leq \pd_{\mI_{n,m-1}}(\bM )
+\pd_{\mI_{n,m}}(\mI_{n,m-1})\leq n+m-1+1=n+m.$$ By (\ref{tordM}),
$$ \pd_{\mI_{n,m}}(M) \leq \max \{ \pd_{\mI_{n,m}}({\rm tor}_\der
(M)), \pd_{\mI_{n,m}}(\bM ) \}\leq n+m,$$ as required. The proof
of the theorem is complete. $\Box $

$\noindent $

{\it Conjecture}. $\gldim (\mI_n)=n$.


\section{The weak and the global dimensions of the Jacobian algebra
$\mA_n$}\label{WGDJJ}

In this section, we prove that the weak  dimension of the Jacobian
algebra $\mA_n$ and of all its prime factor algebras is $n$
(Theorem \ref{J18Oct9}, Corollary \ref{Ja19Oct9}). An analogue of
Hilbert's Syzygy Theorem is established for the Jacobian algebras
$\mA_n$ and for all its prime factor algebras (Theorem
\ref{J19Oct9}).

 A $K$-algebra $R$ has the {\em endomorphism
property over} $K$ if, for each simple $R$-module $M$, $\End_R(M)$
is algebraic over $K$.

\begin{theorem}\label{bigAn}
{\rm \cite{Bav-Ann2001}}  Let $K$ be a field of characteristic
zero.
\begin{enumerate}
\item The algebra $\CA_n$ is a simple, affine, Noetherian domain.
\item The Gelfand-Kirillov dimension $\GK (\CA_n)=3n$ $(\neq
2n=\GK (A_n))$. \item The (left and right) global dimension ${\rm
gl.dim} (\CA_n)=n$. \item The (left and right) Krull dimension
${\rm K.dim} (\CA_n)=n$. \item Let ${\rm d}={\rm gl.dim}$ or ${\rm
d}= {\rm K.dim}$. Let $R$  be a Noetherian $K$-algebra with ${\rm
d}(R)<\infty $ such that $R[t]$, the polynomial ring in a central
 indeterminate, has the endomorphism property over $K$. Then
${\rm d}(\CA_1\t R)= {\rm d}(R)+1$. If, in addition, the field $K$
is algebraically closed and uncountable, and the algebra $R$ is
affine, then
 ${\rm d}(\CA_n\t R)= {\rm d}(R)+n$.
\end{enumerate}
\end{theorem}

$\GK ({\cal A}_1)=3$ is due to A. Joseph \cite{Jos1}, p. 336;
 see also  \cite{KL}, Example 4.11, p. 45.

{\bf The Jacobian algebra $\mA_n$ is a localization of the algebra
$\mI_n$}.  Using the presentations of the algebras $\mI_n$ and
$\mA_n$ as GWAs, it is obvious that the algebra $\mI_n$ is the
two-sided localization, 
\begin{equation}\label{mAnlIn}
\mA_n=S^{-1}\mI_n = \mI_n S^{-1},
\end{equation}
of the algebra $\mI_n$ at the multiplicatively closed subset $S:=
\{ \prod_{i=1}^n(H_i+\alpha_i)_*^{n_i} \, | \, (\alpha_i) \in
\Z^n, (n_i)\in \N^n \}$ of $\mI_n$ where
$(H_i+\alpha_i)_*:=\begin{cases}
H_i+\alpha_i & \text{if }\alpha_i\geq 0,\\
(H_i+\alpha_i)_1& \text{if }\alpha_i<0,\\
\end{cases} $
 since, for all elements $\beta\in \Z^n$,
\begin{equation}\label{vbH}
v_\beta \prod_{i=1}^n(H_i+\alpha_i)_*^{n_i}=\prod_{i=1}^n
(H_i+\alpha_i-\beta_i)_*^{n_i} v_\beta.
\end{equation}
The left (resp. right)  localization of the Jacobian algebra
$$\mA_n = K\langle y_1, \ldots , y_n, H_1^{\pm 1}, \ldots ,
H_n^{\pm 1}, x_1, \ldots , x_n\rangle, \;\; ({\rm where}\;\; y_i:=
H_i^{-1}x_i) $$ at the multiplicatively closed set $S_{y_1, \ldots
, y_n}:=\{ y^\alpha \, | \, \alpha \in \N^n\}$ (resp. $S_{x_1,
\ldots , x_n}:=\{ x^\alpha \, | \, \alpha \in \N^n\}$) is the
algebra 
\begin{equation}\label{mAny}
\CA_n\simeq S_{y_1, \ldots , y_n}^{-1}\mA_n\simeq \mA_n S_{x_1,
\ldots , x_n}^{-1}.
\end{equation}
The algebra $\mA_n$ has the involution $*$. The algebra $\CA_n
\simeq \mA_n / \ga^e$ inherits the involution $*$ since
$(\ga_n^e)*= \ga_n^e$, and so do the algebras $\mA_{n,m}:=
\CA_{n-m}\t \mA_m$ where $m=0, 1, \ldots , n$ and
$\mA_0=\CA_0:=K$. Therefore, the algebras $\mA_{n,m}$ are
self-dual, and so $\lgldim (\mA_{n,m})= \rgldim (\mA_{n,m}):=
\gldim (\mA_{n,m})$.

\begin{theorem}\label{J18Oct9}
Let $\mA_{n,m}:= \CA_{n-m}\t \mA_m$ where $m=0, 1, \ldots , n$ and
$\mA_0 = \CA_0:=K$. Then
 $\wdim (\mA_{n,m})  = n$
for all $m=0, 1, \ldots , n$. In particular, $ \wdim (\mA_n) = n$.
\end{theorem}

{\it Proof}. By Theorem \ref{bigAn}.(1,3) and (\ref{mAny}),
\begin{eqnarray*}
 n &=& \gldim (\CA_n)= \wdim (\CA_n) =\lgldim (S^{-1}_{y_1, \ldots ,
y_n} \mA_n)\leq \wdim  (S^{-1}_{y_1, \ldots ,
y_{n-m}} \mA_n)=\wdim (\mA_{n,m})\\
 &\leq &\wdim (\mI_{n,m})=n\;\;\; ({\rm by}\;\; (\ref{mAnlIn})\;\; {\rm and \;\; Theorem} \;
\ref{18Oct9}).
\end{eqnarray*}
Therefore, $\wdim (\mA_{n,m}) =n$ for all $n$ and $m$.  $\Box $

\begin{corollary}\label{Ja19Oct9}
Let $A$ be a prime factor algebra of the algebra $\mA_n$. Then
$\wdim (A)  =n$.
\end{corollary}

{\it Proof}. By Corollary 3.5, \cite{Bav-Jacalg},  the algebra $A$
is isomorphic to the algebra $\mA_{n,m}$ for some $m$. Now, the
corollary follows from Theorem \ref{J18Oct9}. $\Box $

The next theorem is an analogue of  Hilbert's Syzygy Theorem for
the Jacobian algebras and their prime factor algebras.

\begin{theorem}\label{J19Oct9}
Let $K$ be an algebraically closed uncountable field of
characteristic zero. Let $A$ be a prime factor algebra of $\mA_n$
(for example, $A= \mA_n$) and $B$ be a Noetherian finitely
generated algebra over $K$. Then $\wdim (A\t B) = \wdim (A) +\wdim
(B)  = n+\wdim (B)$.
\end{theorem}

{\it Proof}. Recall that $A\simeq \mA_{n,m}$ for some $m\in \{0,
1, \ldots , n\}$ and $ \wdim (\mA_{n,m}) = n$ (Theorem
\ref{J18Oct9}). Since
\begin{eqnarray*}
n+\wdim (B) &=& \wdim (\mA_{n,m})+\wdim(B) \leq
\wdim (\mA_{n,m}\t B) \\
&\leq &\wdim(S^{-1}\mI_n\t B)\leq  \wdim (\mI_n \t B)\\
& = & n+\wdim (B)\;\; ({\rm by \; Theorem\; \ref{19Oct9}}).
\end{eqnarray*}
Therefore, $\wdim (\mA_{n,m}\t B) =n+\lgldim (B)$. The proof of
the theorem is complete.
 $\Box $

\begin{proposition}\label{e29Nov9}
$n\leq \gldim (\mA_{n,m})\leq n+m$ for all $n\in \N$ and $m=0,1,
\ldots , n$. In particular, $n\leq \gldim (\mA_n) \leq 2n$.
\end{proposition}

{\it Proof}. By Theorem \ref{bigAn}.(3), Proposition
\ref{a29Nov9}, (\ref{mAnlIn}) and  (\ref{mAny}), $n=\gldim
(\CA_n)\leq \gldim (\mA_{n,m}) \leq \gldim (\mI_{n,m})\leq n+m$.
$\Box$

$\noindent $

{\it Conjecture}. $\gldim (\mA_n)=n$.

Department of Pure Mathematics

University of Sheffield

Hicks Building

Sheffield S3 7RH

UK

email: v.bavula@sheffield.ac.uk

\end{document}